\documentclass{amsart}
\usepackage{graphicx}
\usepackage[latin1]{inputenc}
\usepackage[all]{xy}
\usepackage{epsfig}

\input{epsf.sty}

\usepackage{amssymb}
\newtheorem{thm}{Theorem}[section]

\newtheorem{lem}[thm]{Lemma}
\newtheorem{prop}[thm]{Proposition}
\theoremstyle{definition}
\newtheorem{defn}[thm]{Definition}
\theoremstyle{remark}
\newtheorem{rem}[thm]{Remark}
\numberwithin{equation}{section}

\newcommand{\norm}[1]{\left\Vert#1\right\Vert}

\newcommand{\Real}{\mathbb R}

\newcommand{\N}{\mathbb N}
\newcommand{\eps}{\varepsilon}
\newcommand{\To}{\longrightarrow}

\newcommand{\F}{\mathcal{F}}
\newcommand{\id}{{\rm id}}
\newcommand{\Id}{{\rm Id}}
\newcommand{\Homeo}{{\rm Homeo}}
\newcommand{\Lin}{{\rm Lin}}
\newcommand{\Diff}{{\rm Diff}}
\newcommand{\nbhd}{neighborhood }
  
\newcommand{\bilip}{biLipschitz }
\newcommand{\cpt}{compact }
\newcommand{\pt}{point }
\newcommand{\pts}{points }
\begin{document}

\title[Smooth  biLipschitz homogeneous geodesic manifolds]
{Geodesic manifolds with a transitive subset of smooth biLipschitz maps    }%
\author{Enrico Le Donne}%
\address{Yale University,
Mathematics Dept.,
PO Box 208283,
New Haven, CT 06520-8283}%
\email{enrico.ledonne@yale.edu}%


\begin{abstract} 

This paper is connected with the
problem of describing path metric spaces that are homeomorphic
to manifolds and \bilip homogeneous, i.e.,  whose
\bilip homeomorphism group acts transitively.

Our main result is the following.  Let
$X = G/H$  be a homogeneous manifold of a Lie group $G$ and let $d$
be a geodesic distance on $X$ inducing the
same topology.   Suppose  there exists a subgroup $G_S$ of $G$ that
acts transitively on $X$, such that each element $g \in G_S$
induces a locally \bilip homeomorphism of the metric space
$(X,d)$.  Then the metric is  locally \bilip equivalent to a 
sub-Riemannian metric.  Any such metric is defined by a bracket generating $G_S$-invariant sub-bundle of the tangent bundle.

The result is a consequence of a more general fact that requires a transitive family of uniformly biLipschitz diffeomorphisms with a control on their differentials. It will be relevant that the group acting transitively on the space is a Lie group and so it is locally compact, since, in general,  the whole group of \bilip maps, unlikely  the isometry group,  is not locally compact. 

Our method also gives an elementary proof of the following fact. Given a Lipschitz sub-bundle of the tangent bundle of a Finsler manifold, then both the class of piecewise differentiable curves tangent 
to the sub-bundle and  the class  of Lipschitz curves almost everywhere tangent to the sub-bundle 
give rise to the same Finsler-Carnot-Carath\'eodory metric, under the condition that the topologies 
induced by these distances coincide with the manifold topology.

\end{abstract}
\maketitle

\section{Introduction}

\label{Historical context} 
In the last 20 years there has been a surge of interest in the geometry of non-smooth  spaces and in their corresponding \bilip   analysis.   
This movement arose from the interaction between active areas of mathematics concerning the theory of Analysis on Metric Spaces 
 \cite{Semmes, Heinonen-Koskela, Cheeger, Ambrosio-Kirchheim, Lang-Plaut, Heinonenbook, Laakso},   Geometric Analysis  \cite{Gromov-Schoen, Burago-Gromov-Perelman,
Cheeger-Colding}
along with Geometric Group Theory, Rigidity,
and Quasiconformal Homeomorphisms.   
One focus of 
this research has been the study of mappings between non-Riemannian metric structures such as Carnot groups and boundaries of hyperbolic groups
\cite{Pansu, Bestvina-Mess, Margulis-Mostow, Kleiner-Leeb, Gromov, Bourdon-Pajot, Kapovich-Benakli}.

\label{Objectives}
In the present paper, we focus   on the rigidity of  certain non-smooth metric structures on manifolds, namely, geodesic metrics on manifolds that have a transitive group of \bilip homeomorphisms; for short they are called \bilip homogeneous geodesic manifolds.
Every  known example  is locally \bilip equivalent to a homogeneous space $G/H$ equipped with a Carnot-Carath\'eodory metric; here $G$ is a connected Lie group and $H$ is a closed subgroup. Any such metric, also called a sub-Riemannian metric, is defined by a bracket generating sub-bundle of the tangent bundle, also known as completely non-holonomic distribution. For surveys of this area, including how the jargon interchanges between 
sub-Riemannian geometry and  Carnot-Carath\'eodory geometry, see  \cite{Burago, Montgomery} and the papers \cite{Gromov1, Mitchell, Bellaiche}. 

All of the  $2$-dimensional examples, known so far, are locally \bilip equivalent to  the Euclidean plane. 
One may therefore ask whether these are in fact the only examples.
Our main goal in this  paper is to show that under additional assumptions
this is indeed the case.



\label{Motivation} 
The motivation for this question comes  from several sources. First, one can view this as an analogue, in the \bilip category, of Hilbert's fifth problem  on the characterization of Lie groups, solved in \cite{mz}, or the conjectural Bing-Borsuk  characterization of topological manifolds.
Another source of motivation is the work of Berestovski\u\i\; \cite{b, b1, b2}, who showed that a finite dimensional geodesic metric space with transitive  isometry group is isometric to an example $G/H$ as above, except that in the general case one has to use a Finsler-Carnot-Carath\'eodory metric as opposed to a Riemannian-Carnot-Carath\'eodory metric.
In addition, a coarse version of this question in the two dimensional case has arisen in several situations in Geometric Group Theory, and is directly related with the problem considered. For example in \cite{Kapovich-Kleiner}, with the purpose of analyzing
$2$-dimensional Poincar\'e duality groups over commutative rings, 
 Michael Kapovich and Bruce Kleiner studied
quasi-homogeneous quasi-planes, i.e., simply-connected metric cell complexes satisfying a coarse Poincar\'e 
duality  in dimension $2$ on which the group of quasi-isometries acts transitively. Quasi-planes also appeared in  \cite{Kapovich-Kleiner2}.

\label{Considerations} We consider   length metrics since otherwise  there
are many metrics with transitive isometry group that are not
locally \bilip to the standard one, e.g., even on the real line $\Real$, all distances $d(s,t) :=  \sqrt[d]{|s-t|}$, for $d>1$, are translation invariant. In fact,  \bilip homogeneous curves have been studied deeply in \cite{Bishop, GH}.

\label{What we show}
We   show that, locally, the examples mentioned above are the only examples if 
the \bilip maps giving homogeneity come from a Lie group acting by diffeomorphisms.  This assumption is equivalent to  the space  being
homeomorphic to a homogeneous space $X = G/H$, with
  $G$ containing a transitive subgroup of \bilip homeomorphisms. To be precise,  we prove the following:
{\thm \label{thm1} Let $G$ be a Lie group and $H$ be a closed subgroup. Let $X=G/H$ be the corresponding homogeneous manifold equipped with a geodesic distance $d$ inducing its natural 
topology. 
 Suppose there exists a subgroup $G_S$ of $G$ that acts transitively on $X$, and that
acts by maps that are locally \bilip with respect to $d$. 
Then there exists a completely non-holonomic $G_S$-invariant distribution, 
such that any  Carnot-Carath\'eodory    metric coming from it gives a metric that is locally \bilip equivalent to $d$.}
\\

\label{Remarks}
Notice that we do not assume uniform  bounds on the \bilip constants.
Indeed,    if one assumes that $X$ is a geodesic metric space with a transitive
group $G$ of $L$-biLipschitz homeomorphisms, then by taking the supremum over the $G$-orbit
of the distance function, and then the associated path metric, one gets an $L$-biLipschitz equivalent metric
with respect to which $G$ acts by isometries. One can then apply Berestovski\u\i's result in \cite{b2} mentioned above.
However, without  the extra hypothesis about uniformity of \bilip
constants, the argument breaks down altogether. 
Our work can be considered as the first step toward a coarse version of Berestovski\u\i's result. 
In fact, the main steps of our proof and Berestovski\u\i's strategy share some 
common features, although his method is more algebraic.

\label{Assumptions} Our assumption is connected with the fact that, in general, the full group of \bilip maps is not a locally compact group.  
On the other hand, as a consequence of the Ascoli-Arzel\`a Theorem, the isometry group of the space $X$ is a locally compact topological group. 
Berestovski\u\i\; uses this fact in the case that
the isometry group acts transitively to apply the celebrated Gleason-Montgomery-Zippin Theorem and subsequent work \cite[Chapter III]{mz}. The result is a reduction  to the case when the action of the isometry group is topologically conjugate to a transitive smooth action of a Lie group on a smooth manifold, in fact on a homogeneous space $X=G/H$. Thus, in the case of \bilip homogeneous geodesic manifolds, we shall assume that we already have a similar structure, in the sense that the \bilip maps giving homogeneity are coming from a Lie group $G$ acting on a quotient $G/H$.

%

\label{Reformulation}
Since the problem is in fact local, another formulation is as follows. 
Suppose we have a geodesic metric on a \nbhd of the 
origin in $\Real^n$, and a collection of \bilip maps that sends the origin to any point 
in this
neighborhood. Furthermore, suppose that these biLipschitz maps
are, in fact, elements of a smooth ``local'' action.  
Then we can conclude that in a \nbhd of the origin this metric is \bilip equivalent to a 
Carnot-Carath\'eodory metric. 
See the next section for more general statements.

Here is a concrete application of 
 Theorem \ref{thm1}.
\proof[Example:  Affine maps giving homogeneity.] Suppose we have a geodesic distance on the plane such that for any two points there exists an affine map 
that
is locally biLipschitz with respect to the geodesic distance and sends
the first point to the second. 
 Then we can conclude that the distance is locally \bilip equivalent to the Euclidean one.

Indeed, one can apply Theorem \ref{thm1} where $G$ is the group of affine diffeomorphisms and $H = GL_2 (\Real)$. The group $G_S$ is the intersection between $G$ and the group of maps that are locally \bilip with respect to $d$. The metric is in fact Riemannian since there is no proper sub-Riemannian structures in dimension $2$.\\

Several references \cite{b, Gromov1, Montgomery} state that in the definition of Carnot-Carath\'eodory metrics, one gets the same metric when 
considering
either piecewise continuously differentiable, or Lipschitz (or absolutely continuous curves), as horizontal curves.
However, we could not find any proof of this fact in the literature. 
An element of the proof of Theorem \ref{thm1}, detailed in step 5 below, can be 
used to give a simple proof of this fact in the case that the two topologies are assumed to 
coincide.
See the Appendix for details. 

\label{Theorem}
{\thm \label{rem} Let $X$ be a Finsler manifold, equipped with a locally Lipschitz sub-bundle (of the tangent bundle). 
Then both the class of piecewise $C^{1,1}$ curves tangent to the sub-bundle and the class of Lipschitz curves almost everywhere tangent to the sub-bundle
give rise to the same Finsler-Carnot-Carath\'eodory metric, under the condition that the topologies 
induced by these distances coincide with the manifold topology. }

\label{An outline  of the proof}
\subsection{An outline  of the proof of Theorem \ref{thm1}}  $ $
\begin{itemize}
	\item[Step 1.]  Argue that we can assume that   the group $G$ is embedded and closed in the  homeomorphism group of the space $X=G/H$. Thus, every time we will apply the Ascoli-Arzel\`a theorem, the limits of $C^0$-converging subsequences
       will still be   elements of the group and the convergence will be, in fact, $C^\infty$.
 	\item[Step 2.]  Apply a Baire category argument to get a locally transitive set of
       elements of the group that are $C^k$ close to $\id_X$ and are uniformly \bilip with respect to both the metric  $d$ and any Riemannian metric on $X$ that we fixed. 
 	\item[Step 3.]  Prove that, locally,   the distance $d$ is greater than a multiple of some smooth Riemannian distance. 
       Therefore, the geodesics for $d$ are Lipschitz maps for the smooth distance; thus, they are differentiable almost everywhere.
  \item[Step 4.] Define a sub-bundle of the tangent bundle related to the set of velocities of the geodesics. Use it to define a 
       Carnot-Carath\'eodory metric $d_{CC}$. It will be easy to argue that $d \geq  {\rm Constant} \cdot d_{CC}$
       locally.
 	\item[Step 5.] Prove that $d \leq {\rm Constant} \cdot d_{CC}$ locally. 	      
\end{itemize}

\subsection{Organization of the paper} In the next section, we discuss generalizations and variations of Theorem \ref{thm1}., 
Indeed, steps 3, 4, and 5 above constitute a more general fact, specified in 
Theorem \ref{smoothgeneral}, that implies Theorem \ref{thm1}.
In Section 3, we  show  that locally
the distance $d$ is greater than a multiple of any smooth distance. An immediate consequence is
 the absolute continuity of geodesics: any curve that is rectifiable with respect to the geodesic distance $d$ is differentiable almost everywhere.
In Section 4,  we define
a subset of the tangent bundle, 
 and prove that it is, in fact, a sub-bundle; this allows us to define an associated Carnot-Carath\'eodory metric. 
In Section 5, we prove the \bilip equivalence between the   geodesic metric $d$, and the newly defined Carnot-Carath\'eodory metric. 
In    Section 6, we complete steps 1 and 2   above to show how Theorem \ref{thm1}
follows from Theorem \ref{smoothgeneral}.
In Appendix A, we repeat the argument to prove Theorem \ref{rem}. 

\subsection*{Acknowledgements} It is a pleasure to thank the Department of Mathematics of Yale University for the warm and friendly atmosphere that I am enjoying during my Ph.D. program. Above all, I particularly wish to thank Bruce Kleiner for his confidence and support. 

\tableofcontents

\section{The general criterion and some consequences} 
Theorem  \ref{thm1} is a consequence of a more general fact.  
Given a geodesic metric on a domain in $\Real^n$,
whenever there is a transitive family of uniformly biLipschitz diffeomorphisms with a control on their differentials, then the metric is locally biLipschitz equivalent to a Carnot-Carath\'eodory metric. 
Namely, this conclusion can be reached when there is  a family $\F$ of $C^2$-diffeomorphisms and a  base  point $p_0$ such that  the orbit $\F(p_0)$ is a neighborhood of $p_0$, the elements of $\F$  are uniformly \bilip for both the distance $d$ and a fixed Riemannian metric, and the family of the differentials   has a uniform modulus of continuity.



We will denote by $ \norm{A}_{\Lin(\Real^n)}$ the norm as linear operator of a matrix $A\in\Lin(\Real^n)$.

{\thm \label{smoothgeneral} 
Let $N\subset\Real^n$ be a compact  neighborhood of $0\in\Real^n$, equipped with a geodesic distance $d$ that induces the standard topology. 
Suppose there exists a family $\F\subset C^2(N,\Real^n)$ of local $C^2$-diffeomorphisms satisfying:
\begin{description}
\item[Homogeneity] for any $p\in N$, there is $f\in\F$, so that $f(0)=p$; 
\item[Uniform biLipschitz] there exists   $k\in\Real$ such that each $f\in\F$ is a $k$-\bilip map on $N\cap f^{-1}(N)$,   with respect to both the  Euclidean metric and 
the  distance $d$;
\item[Uniform modulus of continuity of derivatives] the family
$\{df\}_{f\in\F}$ has a uniform modulus of continuity, 
i.e., there exists an increasing function $\eta:[0,\infty)\to[0,\infty)$, with $\eta(0)=0$, so that, for any $f\in\F$ and $x,y\in N$,
\begin{equation}\label{equiC1} \norm{(df)_x-(df)_y}_{\Lin(\Real^n)}\leq\eta(|x-y|);\end{equation}
\end{description}
 Then there exists an $\F$-invariant, $C^1$ sub-bundle $\Delta\subset TN$ such that if $d_\Delta$ is any sub-Riemannian metric coming from $\Delta$, then  the geodesic metric $d$ is locally \bilip equivalent to $d_\Delta$.
 }
 
The conclusion of the above theorem is that the regularity of the bundle is $C^1$. However, we will prove that the distribution of Theorem \ref{thm1} is smooth in Proposition \ref{smoothdistr}.\\

The following variation shows how Theorem \ref{smoothgeneral} can be used:

{\thm Let $d$ be  a geodesic metric on a \nbhd $N$ of the 
origin in $\Real^n$. Suppose that there exists a $C^2$ map 
\begin{eqnarray*}
F:\Lambda\times U&\to& N\\
(\lambda,p)&\mapsto &F_\lambda(p),
\end{eqnarray*}
 where $U\subset N\subset\Real^n$ and $\Lambda\subset\Real^N$ are neighborhoods of the origin, with the property
 that, for some neighborhood of the origin $\Lambda_0\subset\Lambda$, the set 
$$\{F_\lambda(0) : F_\lambda\text{ is \bilip w.r.t. } d, \,\lambda\in\Lambda_0\}$$ is a \nbhd of the origin.
Then  in a (possible smaller) \nbhd of the origin the metric $d$ is \bilip equivalent to a 
Carnot-Carath\'eodory metric.}

\proof Clearly, we may assume $F_0=\id_U$, and $U$ and $\Lambda_0$ compact sets. 
We want to show that   the hypotheses of Theorem \ref{smoothgeneral} apply. Indeed, first, since the map $F$ is $C^2$, the condition about  
the family of differentials having a uniform modulus of continuity is immediately verified. Next,
call $m$-\bilip those maps that are $m$-\bilip  with respect to both the Euclidean distance and  the distance $d$, and
consider the sets
$$A_m:=\left\{p\in U\;:\; F_\lambda(0)=p, \;F_\lambda \;m\text{-\bilip}, \,\lambda\in\Lambda_0     \right\},$$
Each $A_m$ is closed. Indeed, take $F_ {\lambda_j}(0)\in A_m\to p$, with $F_ {\lambda_j}$ an $m$-\bilip  map, then by Ascoli-Arzel\`a Theorem, some subsequence converges:  $F_ {\lambda_j}\to F$ uniformly and $F$ is still  an $m$-\bilip  map.
Since here we have assumed  $\Lambda_0$ compact, $\lambda_j\to\lambda_\infty\in\Lambda_0$ up to subsequence, so $F=F_{\lambda_\infty}$ so $p\in A_m$, i.e.,  $A_m$ is closed.
So, by Baire Category Theorem, for some $m\in \N$, the set $\{F_\lambda(0) : F_\lambda \; m-\text{\bilip} \}$ is a \nbhd of a point in $U$ that we can assume to be the origin. 
Theorem \ref{smoothgeneral} can be applied to conclude.
\qed
\\

We now give a generalization of Theorem \ref{thm1} that points out the importance of the assumption that the group acting transitively on the space is a Lie group and so is locally compact. 
Notice that the whole group of \bilip maps, unlikely  the isometry group,  is not in general locally compact. 
\label{Corollary}
According to Montgomery-Zippin's work  \cite{mz},  if a locally 
compact group acts continuously, effectively, and transitively on a manifold, then it is a Lie group, 
Here and in what follows, manifolds are supposed to be connected, however, Lie groups can have infinitely many components. 
Thus, Theorem \ref{thm1} yields the following generalization:

{\thm \label{Corollario} Let $(M,d)$ be a  manifold endowed with a geodesic metric (inducing the same topology). Let $G$ be a  locally \cpt group with a countable base.
Let $G\times M\to M$ be a continuous, effective action of  $G$ on $M$. Suppose there exists a subgroup of $G$ acting 
  transitively on $M$ by \bilip maps (with respect to the metric $d$).
Then $(M,d)$ is locally \bilip equivalent to a homogeneous space equipped with a 
Carnot-Carath\'eodory metric.}

\proof   Following, \cite{mz}, any  locally \cpt group $G$ has the property of having an open subgroup $G'<G$ that is the inverse limit of Lie groups; with the language of Montgomery-Zippin's book, $G'$ has property $A$.

First, we claim that, for any $q\in M$, the orbit of $q$ under $G'$, $G'\cdot q$, is open. This is because, called $H$ the stabilizer of the action,
 the projection $G\to G/H$ is open and the orbit action $G/H\to M$ is an homeomorphism by an standard argument \cite[page 121, Theorem 3.2]{Helgason}.

Now we show that the $G'$-action is still transitive. Indeed, fix a point $p\in M$, suppose by contradiction that $G'\cdot p\neq M$. Hence,
$$ M=\left(G'\cdot p\right)\bigsqcup\left( \bigcup_{q\notin G'\cdot p}G'\cdot q\right)$$
is a disjoint union of two non-empty open
sets of $M$. This contradicts the fact that $M$ is connected.

Thus $G'$ satisfies the hypotheses of Montgomery-Zippin's Theorem \cite[page 243]{mz}, so  $G'$ is a Lie group. So  $G'$ doesn't have small subgroups, thus neither does $G$. By Gleason-Yamabe Theorem, cf. \cite[Chapter III]{mz}, $G$ is a Lie group.
\qed

\section{The given metric is greater than a Riemannian metric}
We now start the proof of Theorem \ref{smoothgeneral}. In particular this section is devoted to showing that curves that are rectifiable with respect to the geodesic distance $d$, are differentiable almost everywhere. To prove this we will show (cf. Proposition \ref{rectgeo}) that the given metric $d$ is locally greater than a Riemannian metric.

We begin with a couple of lemmas. The first asserts that, if we just want $C^1$ maps, we may assume continuity of the differentials when the evaluation at   a base point  goes to the base point.

{\lem \label{continu-differ}
Under the assumptions of Theorem \ref{smoothgeneral}, there exists a family 
$\{f_p\}_{p\in N}\subset C^1(N,\Real^n)$, satisfying $f_p(0)=p$, that is uniformly biLipschitz, satisfies the condition that 
the family 
$\{df_p \}$ has a uniform modulus of continuity, and satisfies: 
\begin{description}
\item[Continuity of $(df_p)_0$ at $0$] the map 
\begin{eqnarray*}
N&\to&\Lin(\Real^n)\\
p&\mapsto &(df_p)_0
\end{eqnarray*}
is continuous at $p=0$. In other words, 
\begin{equation}\label{C0ofdf}   \norm{(df_0)_0-(df_p)_0}_{\Lin(\Real^n)}\to0,\qquad \text{as } p\to0.\end{equation}\end{description}
}
\proof
This is another application of the Baire Category Theorem. 
Set $A_0:=N$. For each $j\in \N$, let $\{V_m^{(j)}\}_{m\in\N}$ be a countable cover of $\Lin(\Real^n)$ by balls of diameter $<j^{-1}$, with respect to some fixed metric inducing the standard topology. 
Inductively define sets $A_j \subseteq \Real^n$ as follows. Let $\F$ be the family from Theorem \ref{smoothgeneral}  and let $$A_{j,m}=\{p\in A_j : \exists f\in\F, f(0)=p, (df)_0\in V^{(j)}_m\}.$$
Then $A_j=\cup_mA_{j,m}$.
By the Baire Category Theorem, there is $m_j\in\N$ so that there exists a ball $A_{j+1}$ of radius $<j^{-1}$ in the closure of $A_{j,m_j}$. Let $q\in \cap_jA_j$. By post-composing with a map that sends $q$ to $0$ we may assume that $q=0$. 

For any $p\neq 0$ in $N$, we have $p\in A_j\setminus A_{j+1}$, for some $j$. Then we define $f_p$ as follows. Since $p\in A_j$, there exists sequences  $\{p_k\}\subseteq A_{j-1}$, $\{f_k\}\subseteq\F$ such that  $p_k\to p$, $f_k(0)=p_k$, and $(df_k)_0\in V^{(j)}_{m_j}$. 
Note that the functions $f_k$  are uniformly Lipschitz and converge at $0$, 
and that the sequence 
$\{(df_k)_0\}$ is equicontinuous and uniformly bounded at $0$. Thus, by 
the Ascoli-Arzel\`a Theorem, 
a subsequence of $\{f_k\}$  converges in the $C^1$ topology. 
Let $f_p$ be one  such limit. 
We can apply this same argument to any sequence $\{f_{p_j}\}$ with $p_j\to 0$ and define $f_0$ as one of the limits of some sequence $\{f_{p_{j_i}}\}$ converging in the $C^1$ topology.
 All of the constructed maps $f_p$'s are $C^1$ and still are uniformly \bilip and their derivatives have uniform modulus of continuity.
Moreover, since $\cap_jV^{(j)}_{m_j}$ is a single point, 
$$(df_p)_0\to (df_0)_0, \quad \text{as } p\to0.$$
\qed

  The next lemma gives  uniform control on the deviation of elements in the family $\F$ from their linear approximations.
{\lem[Uniform distortion control] \label{unifcontrol}
If   a family
$\F$ has the property that the family  $\{df\}_{f\in \F}$ is equicontinuous, so (\ref{equiC1}) holds, then there exists an increasing   function $\omega(t)$ such that, for any  element $f\in \F$ and $y\in N$, 
\begin{equation}
\label{omega}
|(df)_0(y)+f(0) -  f(y)  |\leq  \omega(|{y}|),
\end{equation}
 and $\dfrac{\omega(t)}{t}\to 0$ as $t\to 0$. }

\proof By assumption,  there exists an increasing function $\eta:[0,\infty)\to[0,\infty)$, with $\eta(0)=0$, so that, for any $f\in\F$ and $x,y\in N$,
$$\norm{(df)_x-(df)_y}_{\Lin(\Real^n)}\leq\eta(|x-y|).$$
Let  $\omega(t):=\eta(t)t$. Thus we just need to show  (\ref{omega}) for any $f\in \F$ and $y\in N$. Consider the function $t\in\Real\mapsto f(ty)$. By the Fundamental Theorem of Calculus and the Chain Rule,
$$f(y)-f(0)=\int_0^1\dfrac{d}{dt}f(ty)\,dt=\int_0^1(df)_{ty}\cdot y\,dt.$$
So, 
\begin{eqnarray*}
|(df)_0(y)+f(0) -  f(y)  |&\leq&\int_0^1|(df)_{0}\cdot y-(df)_{ty}\cdot y|\,dt\\
&\leq&\int_0^1\norm{(df)_{0} -(df)_{ty}} |y|\,dt\\
&\leq&\int_0^1\eta( |ty|)|y|\,dt\\
&\leq&\int_0^1\eta( |y|)|y|\,dt\\
&= & \omega(|{y}|).\end{eqnarray*}
\qed

This next proposition is the core of the paper. We get the first important relation between the  geodesic distance $d$ and the  Euclidean metric.

{\prop\label{rectgeo} Let $(N,d)$ satisfy the assumptions of Theorem \ref{smoothgeneral}.
Then some multiple of the distance $d$ is greater than the Euclidean one, i.e., there exists $C>0$ such that
\begin{equation} \label{stima-rettif} \norm{p-q}\leq C \cdot d(p,q),\end{equation}
for any $p,q\in N$.}

\proof 
It suffices to prove (\ref{stima-rettif})
for $p,q$ in a neighborhood $V_0$ of the origin. Indeed, since the family $\F$ acts transitively by $k$-\bilip maps with respect to both metrics, then we would have that any point $p\in N$ has a neighborhood $V$ where
$ \norm{q_1-q_2}\leq C k^2 d(q_1,q_2),$
for $q_1,q_2\in V$. 
Indeed, let $f\in\F$ so that $f(0)=p$. Then, for any $q_1,q_2\in V:=f(V_0)$, we have 
$$ \norm{q_1-q_2}\leq k\norm{f^{-1}(q_1)-f^{-1}(q_2)}
\leq C k d\left(f^{-1}(q_1),f^{-1}(q_2)\right) \leq C k^2 d(q_1,q_2).$$
Now, since both distances are geodesic, we would conclude
$ \norm{p-q}\leq C k^2 d(p,q),$
for all points $p,q\in N$.

 
Suppose (\ref{stima-rettif}) is not true. So there is a sequence of pairs of points $p_n, q_n$ where the ratio of the metrics is smaller and smaller:
$$ \dfrac{d(p_n,q_n)}{\norm{p_n-q_n}}<\dfrac{1}{n} ,$$
for all $n\in\N$. 
We can assume $q_n=0$, since we can move it to $0$ using the transitivity of $k$-\bilip maps. Indeed, let $f_n$ be a $k$-\bilip diffeomorphism such that $f_n(q_n)=0$, for each $n\in \N$. Now,
$$\dfrac{1}{n} >\dfrac{d(p_n,q_n)}{\norm{p_n-q_n}}\geq \dfrac{\dfrac{1}{k}d(f_n(p_n),f_n(q_n))}{k\norm{f_n(p_n)-f_n(q_n)}}=\frac{1}{k^2}\dfrac{d(f_n(p_n),0)}{\norm{f_n(p_n)}}.$$
After replacing $p_n$ with $f_n(p_n)$ and possibly changing indices, we get a sequence of \pts $p_n\in N$ with the property
\begin{equation}
\label{seq}\dfrac{d(p_n,0)}{\norm{p_n}}<\dfrac{1}{n}.
\end{equation}
Note that $p_n \to 0$ as $n\to\infty$. Indeed, $p$ is in the bounded set $N$, so $d(p_n,0)\leq \frac{1}{n}\norm{p_n}\leq  \frac{1}{n}{\rm Diam}(N)\to 0$. Moreover, $d$ gives the usual topology, so $p_n \to 0$.

By Lemma \ref{continu-differ}, there is a neighborhood $U$ of   $0$ such that
\begin{equation*}
 \norm{(df_p)_0-(df_0)_0}_{\Lin(\Real^n)}\leq\dfrac{1}{4},\quad\forall p\in U.
 \end{equation*} 
Thus, since $f_0=\id$, we have, for any vector $v\in\Real^2$,
\begin{equation}\label{unquarto}
 \norm{(df_p)_0(v)-v}\leq
  \norm{(df_p)_0-\Id}_{\Lin(\Real^n)} \norm{v}\leq\dfrac{1}{4}\norm{v}\quad\forall p\in U.
 \end{equation} 
 
All of the balls in the argument shall have the origin as center. Consider an Euclidean ball of radius $R$ contained in the neighborhood $U$. 
Since the topologies induced by the two distances are the same, we can find a $d$-ball of radius $r>0$ inside the Euclidean one. We want to use the fact that the points $p_n$ go to zero, and are different from $0$, to construct a `pseudo' geodesic as a chain of a controlled number of points, whose final point $H(p_n)$ will be outside the $d$ ball of radius $r$. 

Let $\omega$ be the function from Lemma \ref{unifcontrol} that controls the extent to which maps in $\F$ fail to be linear at $0$.
Now, fix $n\in \N$ so large that $\omega(\norm{p_n})<\frac{1}{4}\norm{p_n}.$ Since $p_n$ is now fixed, we will drop the subscript $p=p_n$.
Therefore by estimates (\ref{unquarto}) and (\ref{omega}), we have
\begin{equation}\label{unmezzo}
\norm{p+f(0)-f(p)}\leq \dfrac{1}{2}\norm{p},\quad \forall f\in \F.
\end{equation}

Starting with $f_0=\id$, by recurrence, for $j\in\N$, choose $f_j\in\F$ such that $f_j(0)=f_{j-1}(p)$. We claim that there exists $j\leq  \dfrac{2R}{\norm{p}}$  for which $f_j(p)$ 
is not in the Euclidean ball of radius $R$. In fact, each $f_j(p)$  is quantitatively farther from $0$ than $f_{j-1}(p)$. Indeed, we consider the projection $\pi_p(f_j(p))$ of $f_j(p)$ on the direction of $p$ and prove that 
\begin{equation}\label{stimascalare}
\pi_p(f_j(p)):=\dfrac{\langle f_j(p)| p\rangle}{\norm{p}}\geq \dfrac{j}{2}\norm{p}.
\end{equation}
We prove (\ref{stimascalare}) by induction using the estimates (\ref{unmezzo}) and the fact that
$ \langle u| v\rangle\geq \langle w| v\rangle-\norm{u-w}\cdot \norm{v}$, for all $u,v,w\in\Real^2$.
For $j=1$,
\begin{eqnarray*}
\langle f_j(p)| p\rangle &\geq & \langle p+p| p\rangle-\norm{f_1(p) -p-p}\cdot \norm{p}\\
&\geq & 2\norm{p}^2-\frac{1}{2}\norm{p}^2\geq \frac{1}{2}\norm{p}^2.
\end{eqnarray*}
Assume (\ref{stimascalare})  proved for $j$, then
\begin{eqnarray*}
\langle f_{j+1}(p)| p\rangle &\geq & \langle f_j(p)+p| p\rangle-\norm{f_{j+1}(p)-f_j(p)-p}\cdot \norm{p}\\
&\geq & \langle f_j(p)| p\rangle+\langle p| p\rangle-\norm{f_{j+1}(p)-f_{j+1}(0)-p}\cdot \norm{p}\\
&\geq & \frac{j}{2}\norm{p}^2+\norm{p}^2-\frac{1}{2}\norm{p}^2\geq \frac{j+1}{2}\norm{p}^2.
\end{eqnarray*}  
Therefore, for $j^*:= \dfrac{2R}{\norm{p}}$,  we have
$$\dfrac{\langle f_{j^*}(p)| p\rangle}{\norm{p}}\geq \dfrac{2R}{2\norm{p}}\norm{p}=R.$$
So  $f_{j^*}(p)$ is not in  the Euclidean ball of radius $R$.
Let $H(p)$ be the first point in the sequence $\{f_j(p)\}_{j=0}^{j^*}$
that leaves the ball. Such ball contains the ball of radius $r$ with respect to the distance $d$. Therefore,
\begin{eqnarray*}
r&\leq& d(H(p_n),0)\\
&\leq &\sum_{j=0}^{j^*}d(f_{j}(0),f_{j}(p))\\
&\leq &\sum_{j=0}^{j^*}kd(0,p)\\
&\leq &k(j^*+1)d(0,p)\\
&\leq &k\frac{2R}{\norm{p}}d(0,p)+kd(0,p). 
\end{eqnarray*}

Consider now that $p=p_n$ satisfies (\ref{seq}), so we have
$$0<r\leq k\frac{2R}{\norm{p_n}}d(p_n,0)+kd(0,p_n)\leq \frac{2kR}{n}+kd(0,p_n).$$
But then we have a contradiction since the value $r$, on the left, is constant and greater than $0$ while the right hand side goes to $0$, as $n\to\infty$.\qed



{\rem \label{diffy}
From Proposition \ref{rectgeo} we know that some dilation of the distance $d$ is greater than the Euclidean one. 
So, rescaling the metric $d$ if necessary, we may assume that  $\norm{\cdot}\leq  d $ in $N$. 
From this we can conclude that any $d$-geodesic $\gamma$ is a $1$-Lipschitz map with respect to the Euclidean distance, since
$$|t_1-t_2|=d(\gamma (t_1), \gamma(t_2))\geq \norm{\gamma (t_1)- \gamma(t_2)}.$$
More generally, if $\gamma$ is a $d$-rectifiable curve parametrized by (finite) speed, say smaller than $s$, then
$$s|t_1-t_2|\geq L_d(\gamma (t_1,t_2))\geq d(\gamma (t_1), \gamma(t_2))\geq \norm{\gamma (t_1)- \gamma(t_2)}.$$
In other words, $\gamma$ is an $s$-Lipschitz map with respect to the Euclidean distance, so it is Lipschitz in each coordinate.
At this point we are allowed to use a classical fact in Lipschitz analysis, i.e., Rademacher's Theorem: on $\Real$, any Lipschitz function is differentiable almost everywhere.  Hence any $d$-rectifiable curve is differentiable almost everywhere, in particular, it is rectifiable with respect to the Euclidean distance.
}


\section{The construction of the Carnot-Carath\'eodory metric}
To prove Theorem \ref{thm1} we need to find a sub-bundle 
$\Delta$ of the tangent bundle. As a result  of  Proposition \ref{rectgeo} (cf. Remark \ref{diffy}),  we know that any  $d$-rectifiable curve is differentiable almost everywhere, thus it makes sense to look at the set of velocities of $d$-rectifiable curves. For any $p\in N$, we can now define 
a subset of $T_pN$   as
\begin{equation}\label{Delta}
\Delta_p:=\{\dot\gamma(0)\;:\;\gamma(0)=p, \gamma\; d{\rm -rectifiable\;with\; finite \;speed\; and \;differentiable \;at\;}0\},
\end{equation}
and we let $\Delta=\cup_p\Delta_p\subseteq TN$.

In the next lemma we prove some properties of $\Delta$ such as the fact that it is a sub-bundle, together with some control estimates, needed later, on some family of curves $\Gamma$ representative of $\Delta$.
{\lem[Control on curves $\Gamma$  representing $\Delta$] \label{Gamma}  At any point $p\in N$, the set $\Delta_p$ is a   vector space whose dimension is independent of $p$. The set $\Delta=\cup_p\Delta_p\subset TN $ is  invariant under $\F$. 
Moreover, there exists a special class of curves $\Gamma$, and a constant $S>0$, with the following property:  for any $p\in N$ and any $v\in\Delta_p$, there exists  $\gamma\in\Gamma$ such that $\gamma(0)=p $, $\dot\gamma(0)=v$, and for any $t\in \Real$,
\begin{equation}\label{bounded speed}
{\rm Length}_d(\gamma[0,t])\leq S\norm{v}t.
\end{equation}

Moreover, there is  also an increasing function $\omega_\Gamma:\Real\to\Real$ such that, for any $\gamma\in\Gamma$,
\begin{equation}\label{omega Gamma}
\norm{\gamma(t)-\left(\gamma(0)+ \dot\gamma(0)t\right) }\leq\omega_\Gamma(t\norm{\dot\gamma(0)}),\end{equation}
and $\dfrac{\omega_\Gamma(t)}{t}\to 0$ as $t\to0$. }

\proof Take $w_1,\ldots,w_m\in \Delta_0$ a maximal set of linearly independent vectors coming from paths $\gamma_1,\ldots,\gamma_m$. 
We may assume that the $\gamma_j$'s are parametrized by finite $d$-speed $s_j$, and $\norm{w_j}=\norm{\dot\gamma_j(0)}=1$ and Length$_d(\gamma_j[0,t])\leq s_jt$, for all $j=1,\ldots,m$ and all $t$ where the curves are defined.
Note that for each $\gamma_j$, there is an increasing function $\omega_j$ so that 
$$\norm{\gamma_j(t)-\left(\gamma_j(0)+ \dot\gamma_j(0)t\right) }\leq\omega_j(t).$$
Let $S:=\max_js_j$ and $\omega_\Gamma(t):=\max_j\{\omega_j (t)\}$, then each $\gamma_j$ satisfy (\ref{bounded speed}) and (\ref{omega Gamma}).

To prove that $\Delta_0$ is a vector space we will show that, for any $v\in{\rm Span}(\Delta_0)$, a limit of ``zig-zag'' curves, constructed using the $\gamma_j$'s, is still a rectifiable curve and such limit has tangent at zero equal to $v$. 

For simplicity of exposition, assume $v=\frac{1}{2}(w_1+w_2)$. 
 The zig-zag curves are defined recursively by
\begin{equation}\label{zigzag}
 \sigma_\eps(t):=\left\{ \begin{array}{ccl}
\gamma_1(t)&&0\leq t<\eps \\ \\
 f_{\sigma_\eps(n\eps)}(\gamma_i(t-n\eps)),  &&n\eps\leq t<(n+1)\eps
\end{array} \right.
\end{equation}
where $i=n\; ({\rm mod} \;2)$. 
Each curve $ \sigma_\eps(t)$ is $d$-rectifiable with uniformly bounded speed, i.e., (\ref{bounded speed}) holds. By the Ascoli-Arzel\`a Theorem, the curves $\sigma_\eps$ converge, up to subsequence. Moreover the limit, denoted by $\sigma(t)=\sigma^{(v)}(t)$, is a $d$-rectifiable curve parameterized with same finite speed. In other words, (\ref{bounded speed}) holds for $\sigma^{(v)}$ too.

 \begin{figure}[h]
	{\epsfysize=1.9truein\centerline{\epsfbox{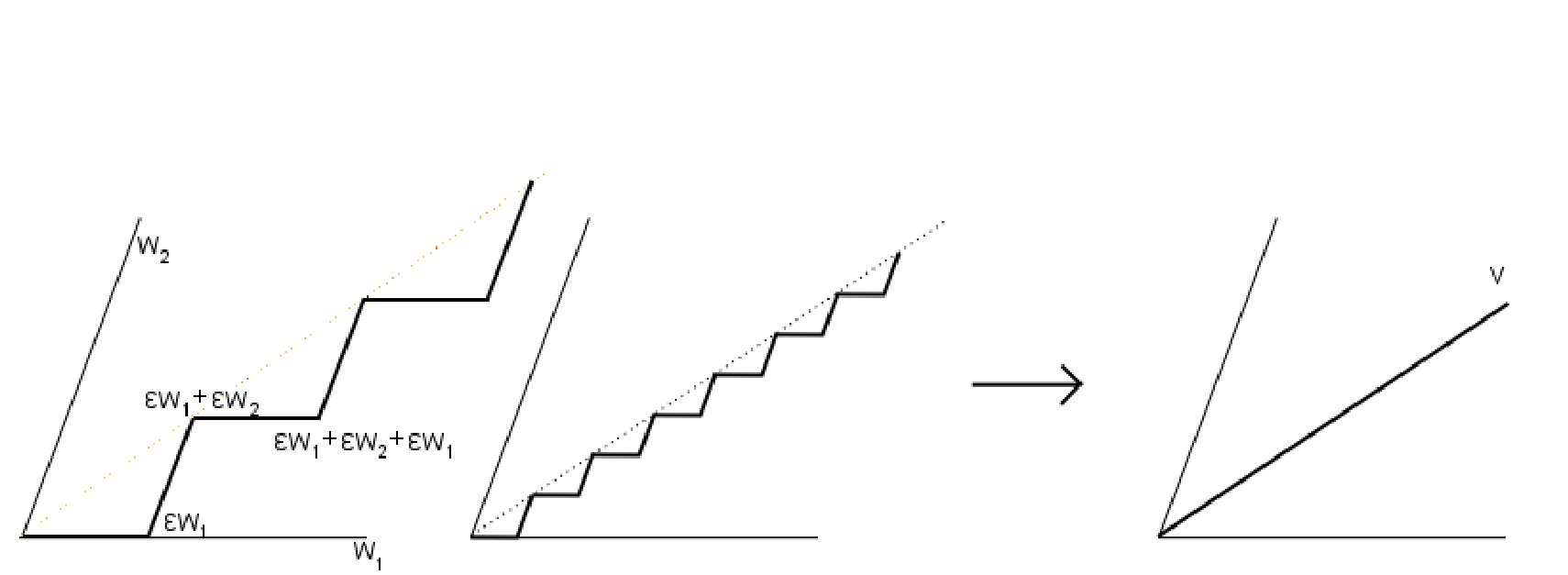}}}
	\caption{The zigzag curves converge to a curve whose tangent at the origin is parallel to the sum of the two vectors $w_1$ and $w_2$.}
\end{figure}
For a better understanding of the curves $ \sigma_\eps(t)$ and their limit, consider first the easier case when 
the $\gamma_j$'s were lines, $\gamma_j(t)=w_jt$, and the $f_p$'s were translations, see Figure 1 above. In this case, the zig-zag curve is
$$\hat\sigma_\eps(t):=\eps w_1+\eps w_2+\eps w_1+\ldots +\bar t w_i,$$
where $t=n\eps+\bar t$ and $i=n\; ({\rm mod} \;2)$. In the sum there are $n+1$ terms. Thus,
$$\hat \sigma_\eps(t)=\dfrac{n}{2} \eps w_1+\dfrac{n}{2} \eps w_2+o(\eps)=\dfrac{t}{2} w_1+\dfrac{t}{2} w_2+o(\eps)=tv+o(\eps)$$ is $\eps$-close to the line $tv$. Therefore, the curve $\hat\sigma(t):=\lim_{\eps\to0}\hat\sigma_\eps(t)$ is such that $$\dfrac{d\hat\sigma}{dt}(0)=v.$$

Define now another auxiliary curve. Let $f_m:=f_{\sigma_\eps(m\eps)}$ and $M_m$ be the matrix $(df_m)_0$. Set
$$\hat{\hat\sigma}_\eps(t):=\eps w_1+\eps M_1w_2+\eps M_2w_1+\ldots +\bar t M_iw_i.$$
By (\ref{C0ofdf}), there is    a neighborhood $U$ of $0$ where $\norm{(df)_0-Id}<\delta$, for all $f\in\F$ with $f(0)\in U$, and so $\norm{(df)_0}<1+\delta$, for all $f\in\F$. Let $t^*>0$ such that 
$$\sigma_\eps(t), \hat\sigma_\eps(t),
\hat{\hat\sigma}_\eps(t)\in U, \qquad\forall t\in [0,t^*].$$
Then we have
\begin{eqnarray*}
\norm{\hat\sigma_\eps(t)-\hat{\hat\sigma}_\eps(t)}&\leq&\sum_{i=1}^n\norm{Id-M_i} \eps \norm{w_i}\\
&\leq&n\eps \delta\\
&\leq&\frac{t}{\eps}\eps\delta\\
&=&t\delta.
\end{eqnarray*}
Also, \begin{eqnarray*}\label{checasino}
\norm{ \sigma_\eps(t)-\hat{\hat\sigma}_\eps(t)}
&=&\norm{  f_n(\gamma_i(\bar t)) - \eps w_1-\eps M_1w_2-\eps M_2w_1-\ldots -\bar t M_iw_i}\\
&\leq& \norm{  f_n(\gamma_i(\bar t))-(d f_n)_0\gamma_i(\bar t)-f_n(0)}+\norm{ (d f_n)_0\gamma_i(\bar t) - \bar t (d f_n)_0\dot\gamma_i(0) } \nonumber\\
&& \qquad \qquad+\norm{ f_{n-1}(\gamma_{i-1}(\bar t)) - \eps w_1-\eps M_1w_-\ldots -\bar t M_{i-1}w_{i-1}} \nonumber\\
&:=& I+II+III.\nonumber
\end{eqnarray*}
Consider the last equation as the sum of three terms: $I$, $II$, and $III$. Note that the first term  can be  bounded, cf. (\ref{omega}), by $$I\leq\omega(\norm{\gamma_i(\bar t)})\leq\omega(\eps).$$ The second term is bounded by
$$II\leq\norm{(d f_n)_0}\norm{\gamma_i(\bar t) - \bar t \dot\gamma_i(0) }\leq(1+\delta)\omega_\gamma(\eps).$$
Finally, the third term is similar to the initial term: the right hand side in the calculation above, 
 except that we have one term less. Therefore, we can iterate the procedure and get
 \begin{eqnarray*}
 \norm{ \sigma_\eps(t)-\hat{\hat\sigma}_\eps(t)}
&\leq& n(\omega(\eps)+(1+\delta)\omega_\gamma(\eps))\\
&\leq&(t+1)\left(\dfrac{\omega(\eps)}{\eps}+(1+\delta)\dfrac{\omega_\gamma(\eps)}{\eps}\right).
\end{eqnarray*}
We are now ready to calculate the derivative at $0$ of $\sigma(t)$. Indeed, for all $t\in U$,
 \begin{eqnarray*}
 \norm{\dfrac{\sigma(t)}{t}-v}&=&\lim_{\eps\to0} \norm{\dfrac{\sigma_\eps(t)}{t}-v}\\
 &\leq&\lim_{\eps\to0} \left(\norm{\dfrac{\hat\sigma_\eps(t)}{t}-v}+ 
 \norm{\dfrac{ \hat\sigma_\eps(t)-\hat{\hat\sigma}_\eps(t) }{t} }+\norm{\dfrac{\hat{\hat\sigma}_\eps(t)-\sigma_\eps(t)}{t}}\right)\\
 &\leq&0+\dfrac{t\delta}{t}+\lim_{\eps\to0} \dfrac{t+1}{t}\left(\dfrac{\omega(\eps)}{\eps}+(1+\delta)\dfrac{\omega_\gamma(\eps)}{\eps}\right)\\
 &=&\delta.
\end{eqnarray*}
So,
for $t\to0$, since $\delta$ is arbitrarily small, we have $\dfrac{d \sigma}{dt}(0)=v.$ 
%
So $v\in\Delta_0$, i.e., $\Delta_0$ is a vector space. 

From the transitive action of \bilip maps we have that
$\Delta$ is  invariant under $\F$, 
so the dimension of $\Delta_p$ is constant. 

%

The class $\Gamma$  is defined to be the curves $\{f_p\circ \sigma^{(v)}\}$ for $p$ in a neighborhood of the origin and $v\in\Delta_0$. Such curves satisfy 
 inequality  (\ref{omega Gamma}) for a suitable $\omega_\Gamma$ since it is true  for the $\gamma_j$'s that are in a finite number, then for the zig-zag limits and finally for all curves in $\Gamma$, using that (\ref{omega}) implies that the $f$'s have a controlled distortion. For  a similar reason, the curves in $\Gamma$ have uniform bound on the speed, i.e., (\ref{bounded speed}) holds.
\qed

From the previous lemma we have that, for each $p\in N$, the set $\Delta_p$ is a $k$-dimensional plane. 
So the map $p\mapsto (p,\Delta_p)$ is a (not necessarily continuous) section of the grassmannian bundle.
In the next subsection we will prove that this map is $C^1$ and so $\Delta$ is  a $C^1$ sub-bundle.

We have noticed, more than once, that the $d$-geodesics in our setting are Lipschitz curves with respect to the Euclidean metric, therefore they are absolutely continuous functions, i.e., they are differentiable almost everywhere and each curve is the integral of its derivative that is a priori just an $L^1$ function. On the other hand, each absolutely continuous curve can be reparametrized to be Lipschitz with respect to the Euclidean metric.

\defn Fixing a distribution $\Delta$, a curve is called \textit{horizontal} if it is absolutely continuous and it derivative lies in the distribution $\Delta$ wherever it exists.

We can now consider another distance on $N$. In the literature, this distance has many different names: \textit{Carnot-Carath\'eodory metric}, \textit{Sub-Riemannian metric, geometric control metric, nonholonomic mechanical metric}.
\begin{equation}\label{d_{CC}}
d_{CC}(p,q)=\inf\left\lbrace {\rm Length}_{\norm{\cdot}}(\sigma)\;|\;\sigma \;{\rm  horizontal \;from}\; p \;{\rm  to} \;q
\right\rbrace,
\end{equation}
where ${\rm Length}_{\norm{\cdot}}$ denotes the length with respect to the Euclidean metric.
%


\subsection{Continuity of the sub-bundle}
To prove that the sub-bundle $\Delta$ is in fact a $C^1$ sub-bundle, we will use a result that give a characterization of  $C^1$ sub-manifold as ambiently $C^1$ -homogeneous compacta.

A set $A\subset \Real^n $
is said to be {\it ambiently $C^1$-homogeneous} if for every pair of points $x , y\in A$, there 
exist neighborhoods $O_x$ and $O_y$ in $\Real^n$ and a $C^1$ diffeomorphism 
$$h : (O_x , O_x \cap A, x )\to(O_y , O_y \cap A, y ).$$ 
{\thm [\cite{Repovs}] \label{compacta} Let $A \subset \Real^n$ 
be compact. Then $A$ is ambiently $C^1$-homogeneous if and only if $A$ is a $C^1 $ submanifold  of $\Real^n$. }

The original proof of this result in  \cite{Repovs} requires the Rademacher Theorem. Shchepin and Repov\v{s} \cite{Shchepin} simplify the proof by eliminating the need to invoke Rademacher. 

Clearly, since the statement of Theorem \ref{compacta} is local, then the assumption of compactness can be replaced by local compactness.

{\prop Let $\Delta\subset TN$ be a   (not necessarily continuous)   distribution. 
Suppose   $\F=\{f_p\}_{p\in N}$ is  a transitive family  of local $C^1$ diffeomorphisms of $N$, $f_p(0)=p$, that leaves invariant $\Delta$, i.e.,
$(df)(\Delta)=\Delta$ for any $f\in\F$.
If the map 
$p\mapsto (df_p)_0$ is continuous at $p=0$, then $\Delta$ is  a continuous sub-bundle. If moreover there is another transitive family of $C^2$ diffeomorphisms  that leaves invariant $\Delta$, then $\Delta$ is
a $C^1$ sub-bundle of the tangent bundle.}

\proof 
The continuity of $\Delta$ at the origin
is consequence of the continuity of 
$p\mapsto (df_p)_0$.
 Then continuity everywhere follows from the invariance under the transitive family of $C^1$ diffeomorphisms.

%
%

Let us denote by $h(p):=\Delta_p$ the $k$-dimensional plane as an element in the $k$-dimensional grassmannian ${\rm Gr}_k(T_pN)$. 
We want to prove that $h$ is a $C^1$ function.
Consider 
the map 
$$\sigma: p\mapsto (p, h(p)),$$
i.e., the induced section of the grassmannian bundle ${\rm Gr}_k(TN)$. 

First notice that, since we proved that $\Delta$ is continuous, then $h$ is  continuous. Being $\sigma$ the graph of the continuous function $h$, then the image $M:=\sigma(N)$ is closed in ${\rm Gr}_k(TN)$. 
Now, if  $f$ is $C^2$, then its differential $df$ induces a $C^1$ map $f_*$ of ${\rm Gr}_k(TN) $.
Moreover if 
a family of $C^2$ diffeomorphisms acts transitively on $N$ preserving $\Delta$, then the induced maps on the grassmannian bundle act transitively on $M$ preserving $M$.
Thus, we   use  Theorem \ref{compacta}: since $M$ is a closed, so locally compact,  subset of the manifold ${\rm Gr}_k(TN)$ that is  ambiently $C^1$-homogeneous, then   $M$ is   a $C^1$ sub-manifold of ${\rm Gr}_k(TN)$. 

Let $V\subseteq N$ be a neighborhood that trivialize the bundle
$${\rm Gr}_k(TV)=V\times\Real^m.$$
Call $\pi:V\times\Real^m\to V$ the projection on the first space. Let $M':=M\cap (V\times \Real^m)$. Since $M$ is   a $C^1$manifold, then 
$$\rho:=\pi_*|_{TM'}:TM' \to TV$$
is a bundle map.
Now, in the case that there exists a point $x\in M'$ such that
\begin{equation}
\label{TMTN}
\rho_x(T_xM')=TN,
\end{equation}
then, by the Implicit Function Theorem, in a neighborhood $U$ of $x$, there is a $C^1$ function $\tilde h:U\to \Real^m$ such that, locally,
$$\left\{(p,\tilde h(p))\;:\:p\in U\right\}=M=\left\{(p,  h(p))\;:\:p\in U\right\}.$$
Thus $h=\tilde h$ and so $h$ is $C^1$ on the neighborhood $U$.
Using again that $M$ is $C^1$-homogeneous, we have that each point in $M$ has a neighborhood on which $h$ is  $C^1$.

If (\ref{TMTN}) is not true for any $x\in M$, we will get to a contradiction.
Indeed, in this case
$$\dim(\ker \rho_x)>0,\qquad\forall x\in M.$$
One then can find an open set $A\subset M$ such that $\dim(\ker \rho_a)=c$, for all $a\in A$, for some constant $c>0$. From general theory of sub-bundles \cite{Atiyah}, we have that 
$\ker \rho\subset TA\subset TM$ is a sub-bundle of $TA\subset TM$ over $A\subset M$.
From this we have that locally there is a non-trivial section $X:A\to \ker \rho$. In other words, $X$ is a vector field on $A$ such that $X_a\in\ker(\pi_*|_{T_aA})$, for all $a\in A$. 
This means that $X_a$ is of the form $(0,X^{(2)}_a)\in TM\times T\Real^m$.
If now $\gamma$ is an integral curve for $X$, then we may assume that $\gamma$ is not constant since $X$ is non-trivial. However,
$$\gamma(t)=\gamma(0)+\int_0^t X_{\gamma(t)}dt =\gamma(0)+\left(0,\int_0^t X^{(2)}_{\gamma(t)}dt \right).$$
But this contradicts the fact that $M$ is a graph.
\qed

\section{Proof of biLipschitz equivalence}
In the previous section we used the fact that the $d$-rectifiable curves are differentiable almost everywhere, by   Proposition \ref{rectgeo}, 
 to construct a distribution $\Delta$ coming from the derivatives of such curves. Now, with the next result, we conclude the proof of Theorem \ref{smoothgeneral}.

{\thm \label{mainstat} Let $(N,d)$ satisfy the assumptions of Theorem \ref{smoothgeneral}. Let $\Delta$ be the distribution defined above in (\ref{Delta}). Let $d_{CC}$ be any Carnot-Carath\'eodory metric associated to  $\Delta$.
Then we have, locally,
$$\frac{1}{L}d_{CC}(p,q)\leq d(p,q)\leq L d_{CC}(p,q),$$
for some $L>0$.}

\proof We need to prove two inequalities.  
\subsection{The first inequality} This is straightforward. Given $p\in N$, let $\gamma_p$ be a $d$-geodesic from $0$ to $p$. 
Since by Proposition \ref{rectgeo} the $d$ distance is greater than the Euclidean distance,  $\gamma_p$ is Lipschitz, thus it is differentiable almost everywhere. We may parametrize $\gamma_p$ by arc length with respect to $d$, so $\gamma_p:[0,T]\to  N$, where $T:=d(0,p)$. We claim that $\norm{\dot \gamma_p}(t)\leq 1$, 
for almost every $t$. Indeed, for any point $t$ of differentiability,
$$\norm{\dot \gamma_p(t)}=\lim_{h\to0}\dfrac{\norm{\gamma_p(t+h)-\gamma_p(t)}}{h}\leq \lim_{h\to0}\dfrac{d(\gamma_p(t+h),\gamma_p(t))}{h}=\dfrac{h}{h}= 1.$$
So
$$d_{CC}(p,0)\leq {\rm Length}_{\norm{\cdot}}(\gamma_p)=\int_0^T \norm{\dot \gamma_p(t)}dt\leq \int_0^T1dt=T=d(0,p).$$

\subsection{The  second inequality}
 
Given a point $p\in N$, we want to construct a $d$-rectifiable curve $\sigma$ that starts at $0$ and ends arbitrarily close to $p$, 
whose $d$-length is close to the $CC$-distance of $p$ from $0$. 
This will be enough since the metric $d$ gives the standard topology.
To construct such a curve, we will use the curves of the family $\Gamma$ defined in Lemma \ref{Gamma}. For any $v\in\Delta$ there is a pre-chosen curve $\gamma_v\in\Gamma$ such that
$\dot\gamma_v(0)=v$, and these curves have a common bound for the speed (\ref{bounded speed}) and for the distance from the linear approximation (\ref{omega Gamma}).

Take any $\eta: [0,T]\to\Real^n$ that is a Lipschitz curve, almost everywhere 
tangent to the distribution $\Delta$, with $\eta(0)=0$, $\eta(T)=p$,  i.e.,  one of the candidate curves in the  calculation of the CC-distance between $0$ and $p$. We can suppose that $\eta$ is parametrized by arc length, i.e.,   $\norm{\dot\eta}=1$, so $T={\rm Length}_{\norm{\cdot}}(\eta)$. Our goal is to show that $T$ is greater than a fixed constant times $d(0,p)$.

Recall that $\Delta$ is $C^1$. Thus, in a neighborhood of $0$, that we will still call $N$, we can find a $C^1$ framing $\{X_1,\ldots,X_k\}$ of $\Delta$, i.e., each $X_j$ is a $C^1$ vector field and
$$\Delta_p={\rm Span}\{(X_1)_p,\ldots,(X_k)_p\},\quad\forall p\in N.$$
We may assume that $\{(X_1)_p,\ldots,(X_k)_p\}$ is an orthonormal basis of $\Delta_p$. Since $N$ is compact, for some $C>0$, each vector field $X_j$ is  $C$-Lipschitz, i.e., 
$$\norm{(X_j)_{q_1}-(X_j)_{q_2}}\leq C \norm{q_1-q_2}, \quad\forall q_1,q_2\in N.$$
A consequence is that if $v\in\Delta_{q_1}$ is such that $\norm{v}\leq 1$ then there is $w\in\Delta_{q_2}$ with $\norm{w}\leq 1$ such that 
$\norm{v-w}\leq C \norm{q_1-q_2}.$
Indeed, since $\{(X_j)_{q_1}\}$ is an orthonormal basis, there are numbers $a_j$ such that 
$v=\sum_ja_j (X_j)_{q_1}$ with $\sum_j |a_j|\leq 1.$ 
Thus $w:=\sum_ja_j (X_j)_{q_2}$ satisfies 
$$|v-w|\leq \sum_j|a_j| \norm{ (X_j)_{q_1}-(X_j)_{q_2}}\leq  C \norm{q_1-q_2}.$$
What we conclude is that each  $v\in\Delta_{q_1}$  with $\norm{v}\leq 1$ has distance less than 
$C \norm{q_1-q_2}$ from the unit ball in $\Delta_{q_2}$. Denoting by $U(\Delta_{q})$  
the unit ball in $\Delta_{q}$, we write
\begin{equation}\label{C-Lip}
{\rm dist} ( U(\Delta_{q_2}),U(\Delta_{q_2}) )\leq C \norm{q_1-q_2}
\end{equation}


\subsubsection{The construction of $\sigma$.}
Take $\eps>0$. Construct piece-by-piece a curve $\sigma$ in the following way. 
Start at $0=\eta(0)$. After a suitable choice of a  vector $v_0\in \Delta_0$, we will take the curve $\gamma_{v_0}(t)\in \Gamma$, where $\Gamma$ is the fixed family of curves from Lemma \ref{Gamma}, and then we will define the first piece of $\sigma(t)$ as, for $0\leq t \leq \eps$,
\begin{equation}\label{parte1}\sigma(t):=\gamma_{v_0}(t).\end{equation}
By (\ref{C-Lip}), since $\norm{\dot\eta(t)}\leq 1$ for a.e.~$t$, we have that, for a.e.~$t\leq\eps$,
$${\rm dist} ( U(\Delta_{0}),\dot\eta(t)  )\leq C \norm{\eta(t)}\leq C|t|\leq C\eps,$$
since $\eta$ is parametrized by arc length.
Since the unit ball $U(\Delta_{0})$ is convex, we have 
$${\rm dist} \left(U(\Delta_0), \frac{1}{\eps}\int_0^\eps\dot\eta(t)dt\right)\leq C\eps.$$

Therefore there exists a  $v_0\in \Delta_0$ with $\norm{v_0}\leq1$ such that \begin{equation}\label{extv0}
\norm{v_0 - \frac{1}{\eps}\int_0^\eps\dot\eta(t)dt}\leq  C\eps.
\end{equation}
So for $0\leq t \leq \eps$, the curve $\sigma$ has been defined by (\ref{parte1}).

For the inductive construction of $\sigma$ suppose that for any $t\leq n\eps$, the value $\sigma(t)$ has been defined. 
We shall define $\sigma$ as, for $n\eps<t\leq (n+1)\eps$,
$$ \sigma(t):=\gamma_{v_n}(t-n\eps),$$
for a suitable choice of $v_n\in \Delta_{\sigma(n\eps)}$ and its related $\gamma_{v_n}\in \Gamma$.

First note that $\lim_{t\to n\eps^+}\sigma(t)=\gamma_{v_n}(n\eps-n\eps)=\gamma_{v_n}(0)=\sigma(n\eps)$. Therefore the new piece agrees with the previous one, i.e., the path is continuous. 
Moreover,
$ \sigma((n+1)\eps)=\gamma_{v_n}(\eps).$
Then calculate the (right)-derivative at $n\eps$:
\begin{equation*}
\left.\frac{d}{dt}\sigma(t)\right|_{t=n\eps^+}	=\dot\gamma_{v_n}(n\eps-n\eps)=\dot\gamma_{v_n}(0)		=v_n.		
\end{equation*}
Again using (\ref{C-Lip}), we have that there exists   $w_n\in \Delta_{\eta(n\eps)} $ with $\norm{w_n}\leq1$ such that \begin{equation*}
\norm{w_n - \frac{1}{\eps}\int_{n\eps}^{(n+1)\eps} \dot\eta(t)dt}\leq  C\eps.
\end{equation*}
Also, since $${\rm dist}(U(\Delta_{\sigma(n\eps)}) , w_n)\leq C\norm{\sigma(n\eps)-\eta(n\eps)},$$
there exists a vector $v_n\in \Delta_{\sigma(n\eps)}$  with $\norm{v_n}\leq1$ such that 
$$
\norm{w_n-v_n }\leq C\norm{\sigma(n\eps)-\eta(n\eps)}.
$$
So
\begin{equation}\label{extvn}	
\norm{v_n- \frac{1}{\eps}\int_{n\eps}^{(n+1)\eps} \dot\eta(t)dt }\leq C\norm{\sigma(n\eps)-\eta(n\eps)}+C\eps.
\end{equation}

Let us now estimate $\norm{\eta(T)-\sigma(T)}$. We will show that 
we have a system of the following type:
\begin{equation}\label{extim}\left\lbrace
\begin{array}{ccl}\norm{\eta(\eps)-\sigma(\eps)} &\leq&o(\eps)\\
\norm{\eta(n\eps)-\sigma(n\eps)} &\leq&(1+C\eps)\norm{\eta((n+1)\eps)-\sigma((n+1)\eps)}+o(\eps),\quad\forall n\in\N,
\end{array}\right.
\end{equation}
where $\dfrac{o(\eps)}{\eps}\to 0$, as $\eps\to 0$.
Observe that a sequence of the form
\begin{equation}\label{sistem}
\left\lbrace 
\begin{array}{cl}
a_1= & \alpha \\ 
a_n =& \beta a_{n-1}+\alpha
\end{array} \right. 
\end{equation}
has solution $a_n = \alpha (\beta^{n-1}+\ldots+1)=\alpha \dfrac{1-\beta^n}{1-\beta}.$
So, from (\ref{extim}),
\begin{eqnarray*}\norm{\eta(T)-\sigma(T)}&=&\norm{\eta\left(\frac{T}{\eps}\eps\right)-\sigma\left(\frac{T}{\eps}\eps\right)}\\
&\leq&o(\eps)\dfrac{1-(1-\eps C)^\frac{T}{\eps}}{1-(1-\eps C)}\\
&\simeq &\dfrac{o(\eps)}{\eps} (1-e^{CT})\To 0,\qquad {\rm as} \;\eps\to0.
\end{eqnarray*}

\subsubsection*{One big triangular inequality}
Now, let us do the calculation showing (\ref{extim}).
The case $n=1$ is shown by considering the following four curves and comparing them at time $t=\eps$:
\begin{enumerate}
\item $\eta(t)$,
\item $t\frac{1}{\eps}\int_0^\eps\dot\eta(s)ds$
\item  $tv_0$, 
\item $\sigma(t)=\gamma_{v_0}(t)$.
\end{enumerate}
Step by step,
\begin{description}
\item[1 and 2] At time $\eps$, the curves are at  the same point, by the Fundamental Theorem of Calculus.
\item[2 and 3] By (\ref{extv0}), we have
$$\norm{\eps v_0 - \eps\frac{1}{\eps}\int_0^\eps\dot\eta(t)dt}\leq  C\eps^2.$$ 
\item[3 and 4] Since $\gamma\in\Gamma$, by (\ref{omega Gamma}) we have $\norm{\eps v_0-\gamma_{v_0}(\eps)}<\omega_\Gamma(\norm{v_0}\eps)<\omega_\Gamma( \eps)$.
\end{description}
 Thus putting everything all together with the triangle inequality:
 $$\norm{\eta(\eps)-\sigma(\eps)}\leq C\eps^2+ \omega_\Gamma( \eps).$$

For $n>1$, more estimates are needed. We compare the following five curves at time $t=(n+1)\eps$:
\begin{enumerate}
\item $\eta(t)$, 
\item $(t-n\eps)\frac{1}{\eps}\int_{n\eps}^{(n+1)\eps} \dot\eta(t)dt +\eta(n\eps)$, 
\item $(t-n\eps)\frac{1}{\eps}\int_{n\eps}^{(n+1)\eps} \dot\eta(t)dt +\sigma(n\eps)$, 
\item $(t-n\eps) v_n+\sigma(n\eps)$,
\item $\sigma(t)=\gamma_{v_n}(t-n\eps)$.
\end{enumerate}
Step by step,
\begin{description}
\item[1 and 2] At time $(n+1)\eps$, as before,  the curves are at  the same point:
$$\eta((n+1)\eps)=\eta(n\eps)+\eps \frac{1}{ \eps} \int_{n\eps}^{(n+1)\eps} \dot\eta(t)dt.$$
\item[2 and 3] One is just a translations of the other by $\norm {\eta(n\eps)-\sigma(n\eps)}.$
\item[3 and 4] As before, by (\ref{extvn}), 
$$\norm{\eps v_n- \eps \frac{1}{\eps}\int_{n\eps}^{(n+1)\eps} \dot\eta(t)dt }\leq C\eps\norm{\sigma(n\eps)-\eta(n\eps)}+C\eps^2.$$
\item[4 and 5] From (\ref{omega Gamma}),  the distance between the fourth and fifth curve is
\begin{eqnarray*}\norm {\eps v_n+\sigma(n\eps) -  \gamma_{v_n}(\eps)  }&=&\norm {\eps \dot\gamma_{v_n}(0)+\gamma_{v_n}(0)  -  \gamma_{v_n}(\eps)  }\\
&\leq & \omega_\Gamma(\norm{v_n}\eps )\\&\leq&\omega_\Gamma( \eps).
\end{eqnarray*}
\end{description}
Thus, putting everything together with the triangle inequality:
\begin{eqnarray*}\norm{\eta((n+1)\eps)-\sigma((n+1)\eps)} &\leq& C\eps^2 +\norm {\eta(n\eps)-\sigma(n\eps)}+ \eps C\norm{\sigma(n\eps)-\eta(n\eps)}+\omega_\Gamma( \eps) \\ 
&\leq& (1+\eps C)\norm {\eta(n\eps)-\sigma(n\eps)}+ C\eps^2  +\omega_\Gamma( \eps). \end{eqnarray*}
Thus, with the terminology of the system (\ref{sistem}), $\beta=1+\eps C$ and $\alpha=C\eps^2 +\omega_\Gamma( \eps)=o(\eps).$ 
Then, as we observed after (\ref{sistem}), $\norm{\eta(T)-\sigma(T)}\to 0$, as $\eps\to0$. This show that we can choose   $\eps$ to have $\sigma(T)$ as close as we want to $\eta(T)$.

Now we calculate $d(0,\sigma(T))$:
\begin{eqnarray*}
d(0,\sigma(T))&\leq&\sum_{n<T/\eps}d(\sigma(n\eps),\sigma((n+1)\eps))\\
&=&\sum d( \gamma_{v_n}(0),  \gamma_{v_n}(\eps))\\
&\leq& \sum_{1}^{T/\eps} S\norm{v_n}\eps\\
&\leq& \sum_{1}^{T/\eps} S \eps=  S \eps \dfrac{T}{\eps}=  S T\\
&=& S  {\rm Length}_{\norm{\cdot}}(\eta),
\end{eqnarray*}
where we used, in order, the triangle inequality, then the definition of $\sigma$, i.e., the fact that $\gamma_{v_n}(t)=\sigma(n\eps+t)$, then that $\gamma_{v_n}$ is $d$-rectifiable parametrized by (uniformly) bounded speed, i.e., (\ref{bounded speed}) holds, then the bound for $\norm{v_n}$.
\qed


\section{The case of \bilip maps  coming from a Lie group action}

We now describe how Theorem \ref{thm1} can be proved using Theorem \ref{smoothgeneral}. What we need to show is that the properties of the transitive action can be improved, i.e., steps 1 and 2 of the outlined argument in the introduction can be done. Let  $G$, $H$, and $d$ be as in Theorem \ref{thm1}.
\subsection{Getting a closed and embedded subgroup of Homeo(G/H)}
Any element of $G$ induces a diffeomorphism of $G/H$. 
Without loss of generality, we can assume that $G$ acts
effectively, so that it may be viewed as a subgroup of
$\Diff(G/H)$:
the space of all $C^\infty$-diffeomorphisms of $G/H$ equipped with the $C^\infty$ topology given by uniform convergence on compact sets of the functions together with all of their derivatives.
So $G$ has two different natural topologies: the first one as a subset of $\Diff(G/H)$ and the second one (weaker) as a subset of $\Homeo(G/H)$:
the space of all homeomorphisms of $G/H$ equipped with the $C^0$-topology, i.e., uniform convergence on compact sets. The first topology is more helpful since it gives control on the derivatives, however, the second one is easier to control by category arguments. 

The following proposition tells us that we may assume that the inclusion $\iota: G\hookrightarrow \Homeo(G/H)$ is an embedding and that $\iota(G)$ is closed. In other words, for any sequence of elements of $G$, viewed as a sequence of maps on $G/H$, that converges uniformly on \cpt sets, the limit map is still an element of $\iota(G)$, and the convergence is, in fact, as elements of $G$, and so the sequence converges as maps in $\Diff(G/H)$.

In general a Lie group $G$ acting on $G/H$ can fail to have the above property. Here is an example. 
Let $G$ be the group of isometries of $\Real^4$ generated
by the translations and a non-closed $1$-parameter subgroup of
$O(4)$.  So $G$ is a connected Lie group of dimension
$5$, acting on $\Real^4$. Thus, to have a group that is closed in $\Homeo(R^4)$, one has to extend the group to a bigger group, in this case the closure of $G$ in Isom$(G)$. What is not trivial in general is that such larger group can be chosen to be still a Lie group.

{\prop \label{closure} Let $G$ be a  Lie group  and $H$ be a closed subgroup.
Then there exists a  Lie group $\hat G$ that extends the action of $G$ on $G/H$ and
is embedded in
$\Homeo(G/H)$ as a closed set.  (Moreover,  $G/H=\hat G/\hat H $, for some closed subgroup $\hat H$.)
}

The rest of this subsection is devoted to the proof of Proposition \ref{closure}. 
Let $X$ be the homogeneous space $G/H$.
 After taking the quotient of  $G$ by the kernel of the action, we can suppose $G$ acts effectively on $X$. 
Then we can replace $G$ by its  universal cover, so it is a simply connected Lie group acting on $X$ effectively in a \nbhd of the identity $e\in G$. 

Let $V$ denote the subspace of vector fields on $X$ that corresponds
to the Lie algebra of $G$.  In other words, for each $\xi \in L(G):=T_eG$,
the one-parameter subgroup of $G$
$$t\longmapsto \exp(t\xi)\in G$$
acts on $X$ by translation. So, for any $x\in X$ and $t\in\Real$, we can consider the flow on $X$
$$\Phi_\xi(t,x):=\exp(t\xi)\cdot x.$$
Differentiating, we obtain a vector field on $X$ that gives the above flow: for $x\in X$,
\begin{equation}\label{elementinV}\xi(x):=d(\exp(t\xi)\cdot x)_{t=0}\in T_xX.\end{equation}
Abusing terminology the vector field is still called $\xi$ since we can identify $V$ and $L(G)$. Indeed, $V$ is isomorphic to  $L(G)$	 as vector spaces (and even the bracket operation, up to sign, is preserved, as shown in  \cite{Helgason}).
In particular, we point out that there is also a one-to-one correspondence  of the above flows with elements in $V$ (or $L(G)$). Indeed, 
\begin{equation}\label{flows}
\Phi_\xi(t,\cdot)=\Phi_{\xi'}(t,\cdot)\in\Homeo(X),\,\forall t\in \Real\quad\Longrightarrow\quad\xi=\xi'\in V,\end{equation}
because of the local effectivity of the action:  for $t$ small enough, $\exp(t\xi)\cdot x=\exp(t\xi')\cdot x\in G/H$  implies $\exp(t\xi)=\exp(t\xi')\in G$ and then $t\xi=t\xi'\in L(G)$ since $\exp$ is a local diffeomorphism at the origin in $L(G)$. 

The vectors in the Lie algebra of $H$ correspond to those vector fields in $V$ that vanish at the origin $[e]\in X$, 
\begin{eqnarray}\label{xi_in_L(H)}
\xi\in L(H)&\Longleftrightarrow&\exp(t\xi)\in H,\quad\forall t\in \Real\nonumber\\
&\Longleftrightarrow&\exp(t\xi)[e]=[e] ,\quad\forall t\in \Real\label{charofH}\\
&\Longleftrightarrow&\Phi_\xi(t,[e])=[e], \quad\forall t\in \Real \nonumber\\
&\Longleftrightarrow&\xi([e])=0. \nonumber
\end{eqnarray}

Note that if  $g \in G$, then the translation
$\tau_g: X \To X$, induced by the left translation, $xH\mapsto gxH$,
preserves the vector fields in $V$; this is just another manifestation
of the adjoint representation\footnote{The map $h\mapsto g^{-1}hg$ is differentiable and fixes the origin. Its differential at the origin is a homomorphism of the Lie algebra called ${\rm Ad}_g$.
The map
\begin{eqnarray*}
G&\To& GL(T_eG)\\
g&\mapsto&{\rm Ad}_g
\end{eqnarray*}
is a representation of $G$ inside the algebra homomorphisms of the Lie algebra.
} of $G$:
we have the formula
$g\exp(\xi)g^{-1}=\exp({\rm Ad}_g\xi)$, see \cite[page 53]{Kna02}, so ${\rm Ad}_g\xi\in L(G)$ is the push-forward vector field.
However, we shall be interested in the fact that $\tau_g$   preserves the flows of vector fields in $V$; indeed, we get $g\exp(t\xi)=\exp(t{\rm Ad}_g\xi)g$  and so
\begin{equation}\label{ad}\tau_g(\Phi_\xi(t,x)))=\Phi_{ {\rm Ad}_g\xi}(t,\tau_g(x)).\end{equation}

The new group $\hat G$ extending the action of $G$ will come from the set of homeomorphisms of $X$, that, as the elements of $G$ in (\ref{ad}), preserve the flows of vector fields in $V$. 

 If $\{g_k\}\subseteq G$ is a sequence that converges as maps in $\Homeo(X)$ uniformly on compact sets to a homeomorphism
$f: X \To X$,
then we claim that $f$ also {\em preserves} $V$, in the sense that for any 
$\xi \in V$
the flow $\Phi_\xi(t,x)$
is conjugated by $f$ to the flow $\Phi_{\xi'}(t,x)$ for some
$\xi' \in V,$ i.e., for any $t\in \Real$, the diagram
$$\begin{array}{cclc}
\Phi_\xi(t,\cdot):&X &\longrightarrow & X \\
& f\downarrow  \quad  & &\quad\downarrow f \\
\Phi_{\xi'}(t,\cdot):& X &\longrightarrow &X
\end{array}$$
 commutes. 
 Since, because of (\ref{ad}), any $g\in G$ preserves $V$, the claim is a consequence of the more general lemma:
 
{\lem \label{C0closed}
The space of homeomorphisms preserving $V$ is $C^0$-closed.}

\proof  Let $\{g_k\}_{k\in N}$ be a sequence of homeomorphisms (not necessarily coming from the $G$ action) preserving $V$, i.e., for any $\xi \in V$ and any $k$ there exists a $\xi_k\in V$ such that the following diagram commutes:
$$\begin{array}{cclc}
\Phi_\xi(t,\cdot):&X &\longrightarrow & X \\
& g_k\downarrow  \quad  & &\quad\downarrow g_k \\
\Phi_{\xi_k}(t,\cdot):& X &\longrightarrow &X
\end{array}.$$
If $g_k\to f$ in Homeo$(X)$, then the above diagram (for fixed $\xi$ and $t$), 
converges uniformly on \cpt sets to 
$$\begin{array}{cclc}
\Phi_\xi(t,\cdot):&X &\longrightarrow & X \\
& f\downarrow  \quad  & &\quad\downarrow f \\
\Phi_{\infty}(t,\cdot):& X &\longrightarrow &X
\end{array}.$$

Recall that $L(G)$ is finite dimensional, 
so, after passing to a subsequence, either $\xi_k$ converges in direction, i.e., the sequence of unit vectors 
$\frac{\xi_k}{||\xi_k||}$ 
converges to some $v_1$ or is zero for all $k$. In the second case, there is nothing else to prove since, for any $k$, $g_k$ is the identity, and so is the limit.
In the first case, we can complete $v_1$ to a basis  $v_1,\ldots,v_n$ of $ L(G)$. 
Defining $ T_k:=\norm{\xi_k}$,
 write $ \xi_k=T_k\sum_{j=1}^na_{kj}v_j$ 
so, for $j=2, \ldots n$, $a_{kj}\To 0$ and $a_{k1}\To 1$, as $k\to \infty$.


Note that $\Phi_\xi(0,\cdot)=\Id$, and so $\Phi_{\infty}(0,\cdot)=\Id$. Pick $p\in X$.
Since $\Phi_{\infty}$ is continuous in $t$, for any $\eps>0$, there exists $\delta$ such that $${\rm Diam}\left( \Phi_{\infty}([0,\delta],p)\right) <\eps.$$
We denoted by Diam the diameter of a set with respect to some fixed metric on $X$ inducing the same topology.
Now, by uniform convergence, for any $\eps>0$, there exists $K\in\N$ such that, for any $k>K$,
$$\Phi_{\xi_k}([0,\delta],p)\subset{\rm Nbhd}_\eps\left( \Phi_{\infty}([0,\delta],p)\right).$$
Therefore, \begin{equation}\label{Ad_{g_k}}
\Phi_{\xi_k}([0,\delta],p)\subset{\rm Nbhd}_{2\eps}(p).
\end{equation}
 
We may assume that $p$ is not a fixed point of the vector field $v_1$. Take a time $t>0$ such that $q:=\Phi_{v_1}(t,p)\neq p$. Suppose we chose $\eps<|p-q|/3$.
Assume by contradiction that $T_k\to\infty$.
Take $k$ big enough such that $\dfrac{t}{T_k}<\delta$.
Since, for $s_k=\dfrac{t}{T_k}$ very small,
\begin{eqnarray*}
\Phi_{\xi_k}(s_k,p)&=&\Phi_{T_k\sum_{j=1}^na_{kj}v_j }(s_k,p)\\
&=&\Phi_{\sum_{j=1}^na_{kj}v_j }(T_ks_k,p)\\
&=&\Phi_{\sum_{j=1}^na_{kj}v_j }(t,p)\To \Phi_{v_1}(t,p)=q.\\
\end{eqnarray*}
But this contradicts (\ref{Ad_{g_k}}), which says that, for all $k$,  the points $\Phi_{\xi_k}(s_k,p)$ lie in a \nbhd of $p$ and so outside the ball $B(q,\frac{\epsilon}{3})$, for how we have chosen $\eps$, and therefore they cannot converge to $q$.

From the contradiction we deduce that the sequence $\xi_k$ is bounded so, after passing to a subsequence, it converges to some $\xi'$ and $\Phi_{\infty}$ has to be the flow of $\xi'$ (by uniqueness of limit). In particular $\xi'$  is uniquely  determined by $\Phi_{\infty}$, by (\ref{flows}).
We proved that every subsequence has a convergent
sub-subsequence, and the limit is independent of the choice of the subsequence; therefore $\xi_k$
actually converges to a fixed $\xi'\in V$, giving the conclusion of the lemma.
\qed

By (\ref{flows}), the vector field $\xi'$ of the lemma is uniquely determined  by $\Phi_{\infty}$ and so by $f$ and $\xi$. Therefore we have a well-defined function $f_*:V\rightarrow V$, such that $$f(\Phi_\xi(t,x)))=\Phi_{f_*\xi}(t,f(x)).$$ Note that this induced map on the space $V$ is functorial, i.e., $f_*\circ g_*=(f\circ g)_*$ for any such maps $f$ and $g$.
If $g$ is an element in $G$, then $g_*={\rm Ad}_{g}$, so $g_*$ is a Lie algebra homomorphisms of $V$.
Now, suppose that $g_k\in\Homeo(X)$ have the similar property  that the maps $\xi\mapsto (g_k)_*\xi$ are Lie algebra homomorphisms of $V$. Then if $g_k\to f$ in $\Homeo(X)$, the map $\xi\mapsto f_*\xi$ is also a Lie algebra homomorphism of $V$,
because  $f_*\xi= \lim_{k\to \infty}(g_k)_*\xi$. In other words, fixing a base for $L(G)$, the maps $ (g_k)_*$ are square matrices converging pointwise to a square matrix $f_*$.

Moreover, if the origin $[e]\in X$ is preserved by $f$, then $f_*$ preserves $L(H)$, i.e., 
\begin{equation}f([e])=[e]\quad \Longrightarrow\quad f_*(L(H))=L(H);\end{equation} the reason is just the characterization (\ref{xi_in_L(H)}): 
$\xi\in L(H)$ if and only if $\xi([e])=0$ if and only if, for every $t\in\Real$, $[e]=f([e])=f(\Phi_\xi(t,[e]))=\Phi_{f_*\xi}(t,f([e]))=\Phi_{f_*\xi}(t,[e])$ if and only if $f_*\xi([e])=0$ if and only if $f_*\xi\in L(H)$.

We can consider the group Homeo$_V$ of homeomorphisms that preserve
$V$, in the sense of the lemma above and induce a Lie algebra homomorphism on   $V$.
{\defn[${\rm Homeo}_V$]\label{HomeoVdef} The set  
${\rm Homeo}_V$ is the group of homeomorphisms $f\in{\rm Homeo}(X)$ such that there exists  a Lie algebra homomorphism $f_*:L(G)\to L(G)$ with the property 
$f(\Phi_\xi(t,x)))=\Phi_{f_*\xi}(t,f(x)) $, or, explicitly, for all $t\in\Real$, $\xi\in V$, and $x\in X$,
\begin{equation}\label{HomeoV}f(\exp(t\xi)x)=\exp(t f_*\xi)f(x).\end{equation}
}

Lemma \ref{C0closed} just says that 
the closure of $G$ in Homeo$(X)$
is contained in  Homeo$_V$, and, more generally, Homeo$_V$ is closed in Homeo$(X)$.

{\lem  The group $\Homeo_V$ is generated by left translations by elements of $G$ and automorphisms
of $G$ that fix $H$. In particular,  any element $f\in\Homeo_V$ can be written uniquely as the composition of a translation and such an automorphism, in fact, if $[g]= f([e])$, with $g\in G$, then 
\begin{equation}f=\tau_{g}\circ\hat \Lambda_{{\rm Ad}^{-1}_g \circ f_*},\label{normalform}\end{equation}
where $\tau_{g}$ is the translation by $g$ and  $\hat\Lambda_{{\rm Ad}^{-1}_g \circ f_*}$ is the map induced on the quotient by the (unique) group automorphism of $G$ with differential  ${\rm Ad}^{-1}_g \circ f_*$.
}

\proof  We first argue that if a map $f\in\Homeo_V$ fixes $[e]$ and $f_*=\Id_V$ then in fact $f=\id_X$. Indeed, we claim that
the set of fixed points $F:=\{ gH\in G/H \,:\,f(gH)=gH\}$ is non-empty,  closed and open, and so it is all of $G/H$, i.e., the function $f$ is the identity. Indeed, $F$ is non-empty since the class of the identity is in it by assumption and it is closed since it is defined by a closed relation. The fact that $F$ is open is a consequence of $\exp$ being locally invertible. Indeed, take any $g'$ that is close enough to $g\in F$ so that it can be written as $g'=\exp(\xi)g$ for some $\xi\in L(G)$. Then 
$$f(g'H)=f(\exp(\xi)gH)=\exp(f_*\xi)f(gH)=\exp(\xi)gH=g'H.$$

It is a classical fact, \cite[page 49]{Kna02}, that since $G$ is simply connected, for any (Lie algebra) homomorphism $\psi$ of $L(G)$, there exists a unique smooth (group) homomorphism $\Lambda_\psi$ of $G$ such that $(d\Lambda_\psi)_e=\psi$. Moreover, in our setting, when $\Lambda_\psi$  is $H$-invariant (so it passes to the quotient $X=G/H$), then $\Lambda_\psi\in\Homeo_V$ and we point out that $(\Lambda_\psi)_*=\psi$. Indeed, since $\Lambda_\psi$  is a homomorphism, 
$$\Lambda_\psi(\exp(t\xi)x)=\Lambda_\psi(\exp(t\xi))\Lambda_\psi(x)=\exp(t \psi\xi)f(x).$$

Suppose now that $f: X \to X$ is any map belonging to $\Homeo_V$. Take $g \in G$ such that that  $[g]= f([e])$ and pre-compose  $f$ by the translation $\tau^{-1}_g$, 
so that $(\tau^{-1}_g\circ f)([e])=[e]$.
Take the automorphism
$\Lambda: G \to G$
whose induced automorphism  $V \to V$ is the inverse
of $(\tau_g^{-1}\circ f)_*$. Explicitly, we take $\Lambda =\Lambda_{((\tau^{-1}_g)_*\circ f_*)^{-1}}=\Lambda_{ f_*^{-1}\circ{\rm Ad}_g}$.  Moreover, $\tau_g^{-1}\circ f$ fixes the origin and so $(\tau_g^{-1}\circ f)_*$ fixes $L(H)$, and so $\Lambda$ fixes $H$. Therefore, passing to the quotient $G/H$, we have a homeomorphism $\hat\Lambda: X\to X$.
Then
$\hat\Lambda \circ \tau_g^{-1}\circ f : X \to X$
fixes $[e] \in X$ and maps each left invariant vector field to itself, i.e., $(\hat\Lambda\circ\tau_g^{-1}\circ f)_*=\Id$.  Hence, from what we showed at the beginning of the proof,
$\hat\Lambda \circ\tau^{-1}_g\circ f$
is the identity, i.e., $f=\tau_{g}\circ\hat\Lambda^{-1}.$
The uniqueness comes from the fact that the intersection between translations and automorphisms  is trivial.
\qed

{\lem  \label{HomeoVLie} The group $\Homeo_V$ is a Lie group of diffeomorphisms and the inclusion $\Homeo_V \hookrightarrow \Homeo(X)$ is an embedding with closed image.}
\proof 
The previous lemma says that every element of $\Homeo_V$ is a diffeomorphism. In fact, the lemma is claiming more: observe that the group of left translations and the group of automorphisms of $G$ fixing $H$ are both Lie groups; the first one is equivalent to $G$ itself, and the second one is a closed subgroup of Aut$(G)$ and Aut$(G)$ is a Lie group since automorphisms of a Lie group come from automorphisms of the Lie algebra, i.e., linear transformations of a finite dimensional vector space. Let $K_1$ be the group of left translations and $K_2$ be the group of automorphisms of $G$ fixing $H$. The previous lemma says that as sets $\Homeo_V=K_1K_2.$ Note that $K_1$ is normal in $\Homeo_V$. Indeed, for any $\tau_g\in K_1$ left translation by $g\in G$ and any $\phi\in K_2$, we have
\begin{eqnarray*}
(\phi^{-1}\circ\tau_g\circ\phi)(xH)&=&(\phi^{-1}\circ\tau_g)(\phi(x)\phi(H))\\
&=&\phi^{-1}(g\phi(x)H)\\
&=&\phi^{-1}(g)(\phi^{-1}\circ\phi)(x)\phi^{-1}(H)\\
&=&\phi^{-1}(g)xH=\tau_{\phi^{-1}(g)}(xH).
\end{eqnarray*}
So $\phi^{-1}\circ\tau_g\circ\phi\in K_1$. Thus, we have $\Homeo_V=K_1\rtimes K_2$ is a semi-direct product of Lie groups, so it is a Lie group.

Now, the fact that the inclusion has closed image is just Lemma \ref{C0closed}. However, we must show that it is an embedding.  This comes from the fact that $V$ is finite dimensional and every sequence of matrices converges $C^\infty$ as soon as it converges point-wise. Indeed, take $f_k\in\Homeo_V$, we need to show that, under the hypothesis that $\{f_k\}_k$ converges in $\Homeo(X)$, then it converges in $\Homeo_V$. (Recall that $\Homeo(X)$ has the $C^0$ topology 
but $\Homeo_V$ has the $C^\infty$ one.) Since  $f_k\in\Homeo_V$, there are associated maps $(f_k)_*\in GL(V)$. The proof of Lemma \ref{C0closed} shows that for any $\xi\in V$ the sequence $(f_k)_*\xi$ converges to $f_*\xi$, where $f\in\Homeo_V$ is the $C^0$-limit of $f_k$.
Since $\{(f_k)_*\}_k$ are linear endomorphisms of the finite dimensional vector space $V$ that converge point-wise, then the convergence is in fact in $C^\infty(V)$. 

So, since by assumption we have $ \xymatrix{f_k\ar[r]^{C^0}& f}$, then in particular we have convergence at the point $[e]=H$, i.e.,
 $ \xymatrix{g_kH:=f_k(H)\ar[r]^{G/H}& f(H)=:gH}$, for some $g_k\in G$.  This means that there exist $h_k\in H$ such that  $ \xymatrix{g_kh_k\ar[r]^{G}& g}$. Call $g':=g_kh_k$, so $[g_k']=[g_k]=f_k([e])$, thus we can use the formula (\ref{normalform}) and have 
 $$f_k=\tau_{g'_k}\circ \Lambda_{{\rm Ad}^{-1}_{g'_k} \circ f_*}.$$
 Now, since $ \xymatrix{g_k'\ar[r]^{G}& g}$, then
 $ \xymatrix{\tau_{g_k'}\ar[r]^{C^\infty}& \tau_g}$ and 
 $ \xymatrix{ {\rm Ad}^{-1}_{g'_k} \ar[r]^{C^\infty}& {\rm Ad}^{-1}_g}$. From this last formula and from the fact that 
 $ \xymatrix{ (f_k)_* \ar[r]^{C^\infty}&f_*}$, we know that, 
defining $\hat f_k:=\tau_{g_k'}^{-1}\circ f_k= \Lambda_{{\rm Ad}^{-1}_{g'_k} \circ f_*}$,
 we have that $ \xymatrix{\hat f_k\ar[r]^{C^\infty}& \tau^{-1}_g\circ f}$, and so that $\{f_k\}_k$ is converging in the $C^\infty$ topology.
\qed

\proof[End of the proof of Proposition \ref{closure}] Consider the space $V$ of vector fields given by (\ref{elementinV}). Let $\hat G:=\Homeo_V$ be the group as in Definition \ref{HomeoVdef}, of the homeomorphisms of the homogeneous space preserving $V$ in the sense of (\ref{HomeoV}) and inducing homomorphisms on the algebra. Then
from Lemma \ref{HomeoVLie}, we know that  $\hat G$
extends the action of $G$, consists of $C^\infty$ diffeomorphisms,  is a Lie group and is a closed, embedded subgroup of $\Homeo(X)$.  
\qed

\rem \label{C^infty convergence} 
By Proposition \ref{closure}, without loss of generality we may assume that $G$ has the following property:
as soon as a sequence of elements $g_n\in G$ converges as maps of $X$ in the $C^0$-topology, then the limit is a map coming from $G$, the convergence is also in the topology of $G$ itself, and, moreover, since the action is smooth, the sequence also converges  in the $C^\infty$-topology.

\subsection{Getting uniformly \bilip maps close to the identity }

Remember that $G$ acts by diffeomorphisms on $G/H$. Suppose now that, for any two points in a \nbhd of the identity, there is an element of $G$, sending the first point   to the second one, that is \bilip with respect to  the geodesic metric $d$. The following lemma tells us that, by staying enough close to the origin $[e]$ , we may assume that these elements of $G$ give maps that are uniformly \bilip and are $C^\infty$-close to the identity as much as we want. Such argument is based on the Baire Category Theorem and has been used several times in the theory of homogenous compacta, e.g.~in  \cite[Theorem 3.1]{MacManus} or   \cite[Theorem 6.1]{Hohti}.

{\lem \label{smooth approx}
 Let $d$ be a   metric on a homogeneous space $X=G/H$ inducing the usual topology. 
Suppose $G$ is closed and embedded in $\Homeo(G/H)$. Let $U$ be a \cpt \nbhd of the origin $[e]\in G/H$ where the elements of $G$ that are \bilip with respect to 
the  metric $d$, act transitively.
Then there exists a constant $k$ such that, for any $\eps>0$, there exists a smaller \nbhd $U_\eps$ of the origin where the elements of $G$ that are $k$-\bilip with respect to both 
a fixed Riemannian metric on $X$ and the  metric $d$
and are $\eps$ close to the identity with respect to a fixed $C^\infty$ distance, act  transitively.  }

\proof
The group $G$ can be seen as a subset of Diff$(G/H)$. 
Since the action is smooth, the topology of $G$ is the same as that of any of those induced by any $C^\infty(G/H)$-distance\footnote{We recall that $f_n$ converges to $f$ in $C^\infty(X)$ if, for any $p\in X$, for any charts $\phi$ and  $\psi$ at $p$ and $f(p)$ respectively, for any \cpt $K\subset \Real^n$ inside the domain of $\phi$, and any multi-index $\alpha\in \N^n$, we have that the associated seminorm goes to zero, i.e., 
$$\sup_KD^\alpha\left( \phi^{-1}\circ f_n\circ \phi -\phi^{-1}\circ f\circ \phi \right) \To 0.$$
     In general, if $( \rho_n)_{n\in N}$ is a sequence of seminorms defining a (locally convex) topological vector space $E$, then
    $$d(x,y)=\sum_{n=1}^\infty \frac{1}{2^n} \frac{\rho_n(x-y)}{1+\rho_n(x-y)}$$
    is a metric defining the same topology.
}. Fix one such distance and fix $\eps>0$.
Consider a cover of $G$ by a countable number of closed balls $D_n$ for $n\in\N$. 
We will say later how to choose the cover depending on $\eps$.
So $G=\bigcup_{n\in\N}D_n$ with $D_n$ closed. 
Since any element of $G$ is smooth, then in particular it is locally \bilip with respect to any Riemannian metric on $X$. Fix one such Riemannian metric. We will use the term $m$-\bilip to describe to maps that are $m$-\bilip for both the Riemannian metric and 
the  metric $d$, and the set of such maps will be denoted as ${\rm biLip}^{m}$.
Now consider the orbit set of $0$ in $U$ under $m$-\bilip maps in $D_n$, i.e.,
\begin{eqnarray*}A_{n,m}&:=&(D_n\cap {\rm biLip}^m )0\\
&:=&\left\lbrace p\in U \;|\;f(0)=p \;{\rm \;for \;some }\;m{\rm -\bilip \;map} \;f\; \in D_n \right\rbrace.\end{eqnarray*}

By transitivity, we have $U=\bigcup_{n\in\N}A_{n,m}.$
We claim that the $A_{n,m}$ are closed. Indeed, take $p_j\in A_{n,m}$ converging to $p\in U$. 
Choose the $f_j\in D_n$ such that $f_j(0)=p_j.$
The $f_j$'s are $m$-\bilip and $f_j(0)=p_j$ converges.
The Ascoli-Arzel\`a argument implies that, after passing to a subsequence, $f_j$ converge to some $f$, uniformly on \cpt sets and the limit function is $m$-\bilip (for both metrics). 
Since $G$ is assumed to be closed in Homeo$(G/H)$, then in fact $f\in G$.
From Remark \ref{C^infty convergence}, the convergence
is in the $C^\infty$ topology, so, since $D_n$ is $C^\infty$-closed, the limit function belongs to it and its value at $0$ is $p$. Therefore, $f(0)=p$ for an $f\in D_n\cap {\rm biLip}^m$. In other words, $p\in A_{n,m}$, so  $A_{n,m}$ is closed.

The Baire Category Theorem implies that one of the $A_{n,m}$
 has non-empty interior. So there exists a \cpt $V\subset U$, which is a \nbhd of some point $q$, such that $V\subset A_{\bar n,\bar m}$ for some $\bar n, \bar m\in\N$.
Consider $f_q\in D_{\bar n}\cap {\rm biLip}^{\bar m}$  such that $f_q(0)=q$.

We claim that we can take $U_\eps:=f_q^{-1}(V)$ as our new \nbhd and 
$k:=\bar m^4$ as our required constant. Indeed, for any \pts $p_1, p_2\in f_q^{-1}(V)$,
for $i=1,2$, $f_q(p_i)\in V\subset A_{\bar n,\bar m}$, so there exists $f_i\in D_{\bar n}\cap {\rm biLip}^{\bar m}$
such that $f_i(0)=f_q(p_i)$. Thus $$p_2=\left( f_q^{-1}\circ f_2\circ f^{-1}_1\circ f_q\right) (p_1)$$
and $f_q^{-1}\circ f_2\circ f^{-1}_1\circ f_q\in G$ is $k$-biLipschitz.
Moreover,
$f_q^{-1}\circ f_2\circ f^{-1}_1\circ f_q$ is $\eps$ close to the identity in Diff$(G/H)$, if we had previously made a good choice of the cover $\{D_n\}$, considering that the function 
$h\circ f\circ g^{-1}\circ h^{-1}$ is continuous in $f,g,h\in G$. 

Now we explain how to choose the cover, given $\eps$. Consider the map 
$$\begin{array}{ccc}C^\infty(X)\times C^\infty(X)&\To & C^\infty(X) \\
(g,h)&\longmapsto &g\circ h\circ g^{-1}.\end{array}$$
It is continuous and sends $C^\infty(X)\times \{\id\}$ to the identity function.  
Given a fixed $\eps$,   there exists a \nbhd $V_1$ of the identity in $C^\infty(X)$ such that $C^\infty(X)\times V_1$ goes into the $\eps$-\nbhd of the identity function (the convergence is in $C^\infty(X)$ with respect to the $C^\infty$-metric that we fixed).
Now consider the   map 
$$\begin{array}{ccc}C^\infty(X)\times C^\infty(X)&\To & C^\infty(X) \\
 (g,h)&\longmapsto& g\circ h^{-1}.\end{array}$$
It is continuous and sends the diagonal $\Delta$ to the identity function.
Given the  \nbhd $V_1$ of before,   there exists a   \nbhd $V_2$ of $\Delta$ that is sent by the map into $V_1$.

So if we had chosen the cover  so that
$D_n\times D_n\subset V_2$ for any $n\in\N$, then,
 for any $f,g \in D_n$, $f\circ g^{-1}\in V_1$. Thus for any $h\in D_n$, we have that $h\circ f\circ g^{-1}\circ h^{-1}$ lies in an $\eps$-\nbhd of the identity. This was what was left to prove.
\qed

\subsection{The end of the proof of Theorem \ref{thm1} using Theorem \ref{smoothgeneral}}
We can focus on a neighborhood $N$ of a point in the manifold that, after fixing a coordinate chart, is a neighborhood of $0$ in  $\Real^n$. We can also transfer the geodesic metric imposing that the chart is an isometry. From now on we will identify the neighborhood in the manifold and that one in $\Real^n$.

We need to construct now a family of maps.
The idea is to use Lemma \ref{smooth approx} to select, for each $p$ in the neighborhood, a \bilip diffeomorphism $f_p$ whose differential  $d(f_p)$  differs from the identity by
an error that depends only on $d(0,p)$,
and that tends to zero, as $d(0,p)$ tends to zero.

Considerthe  neighborhoods $U_{1/n}$ given by Lemma \ref{smooth approx}, we can suppose that $N=U_1$. Note that $\bigcap_n U_{1/n}=\{0\}$.
For $p=0$, choose $f_0=\Id$; for  $p\in U_{1/n}\setminus U_{1/(n+1)}$, Lemma \ref{smooth approx} gives the existence of  a $k$-\bilip map $f_p$ so that $f_p(0)=p$ and which is $1/n$ close to the identity. 

We need to show that the hypotheses of Theorem \ref{smoothgeneral} apply. Since uniformly \bilip homogeneity is clear, we have left to show that condition 
(\ref{equiC1}) holds. 
 It  is satisfied, since in $N$ the second derivatives of the $f$'s are equibounded, say by $C$, so
\begin{eqnarray*}
\left| \dfrac{\partial f_i}{\partial x_j}(x)-\dfrac{\partial f_i}{\partial x_j}(y) \right| &\leq&\int_0^{|x-y|}\dfrac{d}{dt} \dfrac{\partial f_i}{\partial x_j} \left(x+t\frac{x-y}{|x-y|}\right)\,dt\\
&\leq&\int_0^{|x-y|}\sum_{k=1}^n\dfrac{\partial^2 f_i}{\partial x_k\partial x_j}\left(x+t\frac{x-y}{|x-y|}\right)\frac{(x-y)_k}{|x-y|}\,dt\\
&\leq&nC|x-y|.
\end{eqnarray*}

Hence all of the hypotheses of Theorem \ref{smoothgeneral} apply. Thus, there exists a  sub-bundle $\Delta\subset TN$, defined explicitly in (\ref{Delta}), such that if $d_\Delta$ is any sub-Riemannian metric coming from $\Delta$, then  the geodesic metric $d$ is locally \bilip equivalent to $d_\Delta$. The fact that $\Delta$ is invariant under the action of a transitive subset of a Lie group, will imply, by next proposition, that $\Delta$ is not just $C^1$, but in fact smooth. \qed

{\prop \label{smoothdistr} Assume that $S$ is a set of elements of a
Lie group $G$ that is transitive on the space $X:=G/H$. Let  $\Delta$ be a 
 distribution on $X$ preserved by the action of $S$. Then  $\Delta$  must be real analytic.}

 \proof 
Since $S\subset G$ preserves the distribution $\Delta$, then any product of its elements does. Call $G_S$ the group generated by $S$. 
Consider $G_1:=\bar G_S$ the closure of $G_S$ in $G$. The set  $G_1$ is a Lie group and   it preserves the distribution $\Delta$ too. 

Fix a \pt $p\in X$. Look at the orbit map of $p$ under $G_1$, 
$$\Phi: G_1\To X$$
$$g\mapsto g(p).$$
This is a smooth map, so we can take the derivative at the identity, i.e., the differential from the tangent space at the identity $L:=T_eG_1$, the Lie algebra, and the tangent at $p$, $T_pX$,
$$d\Phi_e : L\to T_pX.$$
Since this map is surjective, we can find a subspace $W$ of the Lie algebra of $G_1$, such that the orbit map
restricted to this subspace is an
isomorphism of vector spaces,
$$d\Phi_e|_W : W\to T_pX.$$
Moreover, $\exp(W)$ is, locally, an analytic sub-variety. The orbit map $\Phi$ restricted to $\exp(W)$
gives an analytic map 
$$\Phi|_{\exp(W)}:\exp(W)\to X.$$
By the Implicit Function Theorem, this is locally an analytic isomorphism.

From the hypothesis we know that, for any $g\in G_1$, $\Delta_{g(p)}=g\Delta_p.$
Now,  $g=\exp(v)$, for some $v\in W$, so we have the formula: 
$$\Delta_{(\exp(v))(p)}=\exp(v)\Delta_p.$$
For any other \pt $q\in X$, let $\Phi_{loc}^{-1}$ be a local inverse of $\Phi$ in a \nbhd of $q$, so
$q=\Phi(\exp(v))=(\exp(v))(p)$ and $\exp(v)=\Phi_{loc}^{-1}(q)$. Then
$$\Delta_q=\Delta_{(\exp(v))(p)}=\exp(v)\Delta_p=\Phi_{loc}^{-1}(q)\Delta_p.$$
This implies that $\Delta$ is smooth because $\Delta_q$ depends analytically on $q$.
\qed


\appendix
\section{Equality of piecewise $C^{1,1}$ and Lipschitz Carnot-Carath\'eodory metrics}
One can define the Carnot-Carath\'eodory distance using as horizontal curves either Lipschitz curves tangent almost everywhere or piecewise  $C^{1,1}$ curves, i.e., $C^1$ curves with Lipschitz derivative, tangent to the distribution. 
We prove now that piecewise $C^{1,1}$ horizontal curves or Lipschitz horizontal curves yield to the same Carnot-Carath\'eodory distances, when they both induce the manifold topology.

{\thm Let $X$ be a Finsler manifold and let $\Delta$ be a locally Lipschitz sub-bundle of the tangent bundle. Let $d_{CC}^{ C^{1,1}}$ and $d_{CC}^{Lip}$ be the Finsler-Carnot-Carath\'eodory metrics where the horizontal curves are chosen to be the class of curves that are, respectively, piecewise $C^{1,1}$ and Lipschitz, tangent to the distribution almost everywhere. Suppose that the two topologies induced by the two distances are the same as the topology of the manifold. Then the distances are the same, i.e., $d_{CC}^{C^{1,1}}= d_{CC}^{Lip}$.}

\textit{Proof.}
The fact that $d_{CC}^{C^{1,1}}\geq d_{CC}^{Lip}$ is obvious since all piecewise $C^{1,1}$ curves are Lipschitz. A priori, the infimum over all of the Lipschitz curve can be strictly smaller.

Since both metrics are path metrics, it suffices to prove the statement locally. So we may suppose that we are in $\Real^n$ with a fixed norm $\norm{\cdot}$. In general, the space could be covered by small balls in which charts give a $(1+\eps)$-\bilip approximations. The $\eps$ goes to zero as the diameters of the balls go to zero.

Take $\eta: [0,T]\to\Real^n$ to be a Lipschitz curve, almost everywhere 
tangent to the distribution $\Delta$, with $\eta(0)=0$, $\eta(T)=p$, i.e., one of the candidate to calculate the Lipschitz CC-distance between $0$ and $p$. We can suppose that $\eta$ is parametrized by arclength, i.e.,   $\norm{\dot\eta}=1$ a.e., so $T={\rm Length}(\eta)$.

We will construct a sequence of piecewise $C^{1,1}$ curves whose length is smaller than or equal to the length of $\eta$, going from $0$ to a sequence of points that converges to $p$. Since the topologies are the same, this will give the conclusion.

Take $\eps>0$. Construct piece-by-piece a curve $\sigma$ in a way similar to that of Section 5.2. 
Start at $0=\eta(0)$. After a suitable choice of a  vector $v_0\in \Delta_0$, we will take a curve $\gamma_{v_0}(t)$, 
as the next Lemma \ref{vector field} says, and then we will define the first piece of $\sigma(t)$ as, for $0\leq t \leq \eps$,
$$\sigma(t):=\gamma_{v_0}(t).$$
Each curve $\sigma$ that we will construct  will have length less than $2T$. On the Euclidean ball of radius $2T$, the distribution $\Delta$ is  $C$-Lipschitz, for some $C>0$. By the same discussion that in Section 5.2 led us to (\ref{C-Lip}), we have that there exists
$v_0\in \Delta_0$  with $\norm{v_0}\leq1$ such that \begin{equation}\label{extv02}
\norm{v_0 - \frac{1}{\eps}\int_0^\eps\dot\eta(t)dt}\leq  C\eps.
\end{equation}

For the inductive construction of $\sigma$ suppose that for any $t\leq n\eps$, $\sigma(t)$ has been defined. 
We shall define $\sigma$ as, for $n\eps<t\leq (n+1)\eps$,
$$ \sigma(t):=\gamma_{v_n}(t-n\eps),$$
for a suitable choice of $v_n\in \Delta_{\sigma(n\eps)}$ and its related $\gamma_{v_n}(t)$ given by Lemma \ref{vector field}.


Take  $w_n\in \Delta_{\eta(n\eps)} $ with $\norm{w_n}\leq1$ such that 
\begin{equation*}
\norm{w_n - \frac{1}{\eps}\int_{n\eps}^{(n+1)\eps} \dot\eta(t)dt}\leq  C\eps.
\end{equation*}

Since $\Delta$ is $C$-Lipschitz, there exists 
a vector $v_n\in \Delta_{\sigma(n\eps)}$ with $\norm{v_n}\leq1$ such that 
$$
\norm{w_n-v_n }\leq C\norm{\sigma(n\eps)-\eta(n\eps)}.
$$
So
\begin{equation}\label{extvn2}	
\norm{v_n- \frac{1}{\eps}\int_{n\eps}^{(n+1)\eps} \dot\eta(t)dt }\leq C\norm{\sigma(n\eps)-\eta(n\eps)}+C\eps.
\end{equation}

To  estimate $\norm{\eta(T)-\sigma(T)}$, one uses system (\ref{extim}).
The strategy of the proof for bounding the terms in system (\ref{extim}) is the same as that one in Section 5.2, except that instead of using (\ref{omega Gamma}), one uses   Lemma  \ref{vector field}. Also, estimates  (\ref{extv0}) and (\ref{extvn}) are replaced by (\ref{extv02}) and (\ref{extvn2}).

The conclusion is that, as we observed after (\ref{sistem}), $\norm{\eta(T)-\sigma(T)}\to 0$, as $\eps\to0$. This shows that we can choose   $\eps$ to have $\sigma(T)$ as close as we want to $\eta(T)$.

Now we calculate the length of $\sigma(0,T)$:
\begin{eqnarray*}
{\rm Length}(\sigma[0,T])&\leq&\sum_{n<T/\eps}{\rm Length}(\sigma[n\eps),\sigma((n+1)\eps])\\
&=&\sum {\rm Length}( \gamma_{v_n}[0,\eps])\\
&\leq& \sum_{1}^{T/\eps} \norm{v_n}\eps\\
&\leq& \sum_{1}^{T/\eps} \eps=  \eps \dfrac{T}{\eps}=  T\\
&=& {\rm Length} (\eta[0,T]),
\end{eqnarray*}
where we used, in order, the triangle inequality, then the definition of $\sigma$, then that $\gamma_{v_n}$ is   parametrized by (uniformly) bounded speed, as Lemma \ref{vector field} says. 
\qed

{\lem \label{vector field} Given $v\in \Delta$ there exists a $C^{1,1}$ curve $\gamma_v$ starting at $v$ tangent to the distribution $\Delta$, parametrized by speed smaller than $\norm{v}$ and such that $$\norm{\gamma_v(t)-[\gamma_v(0)+\dot \gamma_v(0)t]}\leq Ct^2. $$}

\proof Extend $v$ to a  vector filed $X\subset\Delta$ using the orthogonal projection:
$$X_p:=\pi_{\Delta_p}(v). $$
The sub-bundle $\Delta$ is $C$-Lipschitz, thus $X$ has the properties of being $C$-Lipschitz, $v\in X$ and $\norm{X_p}\leq\norm{v}$. 
Let $\gamma_v$ be the integral curve of $X$ starting at $v$, i.e., $\dot \gamma_v(0)=v$ and $\dot \gamma_v(t)=X_{\gamma_v(t)}$. Then
\begin{eqnarray*}
{\rm Length}(\gamma_v[0,t])&=& \int_0^t\dot \gamma_v(s)ds\\
&\leq&\int_0^t\norm{\dot \gamma_v(s)}ds\\
&=&\int_0^t\norm{X_{\gamma_v(t)}}ds\\
&\leq&\int_0^t\norm{v}ds\;=\;t\norm{v}.
\end{eqnarray*}
In other words, $\gamma_v$ is of speed smaller than $\norm{v}$. 
The rest of the conclusion of the lemma is clear:  since $X$ is $C$-Lipschitz then it is differentiable a.e.~with derivative bounded by $C$. Thus $\gamma_v$ is a $C^2$ curve a.e.~with second derivative bounded by $C$.
\qed

\nocite{Heinonenbook}
\bibliography{bilip_homog}
\bibliographystyle{amsalpha}
\end{document}